\documentclass[a4paper,twoside,12pt,english]{article}

\usepackage[applemac]{inputenc}
\usepackage[english]{babel}
\makeatletter
\renewcommand\section{\@startsection
{section}{1}{0mm}%
{-2\bigskipamount}%
{\bigskipamount}%
{\normalfont\normalsize\bfseries}%
}

\newcommand\dateymd{\number\year, \ifcase\month\or
January\or February\or March\or April\or May\or June\or
July\or August\or September\or October\or November\or
December\fi, \number\day}
\newcommand\printtime{%
\c@hours=\time \divide\c@hours by60
\c@minutes=\c@hours \multiply\c@minutes by-60
\advance \c@minutes by \time
\ifnum\c@hours<10 0\fi\the\c@hours:%
\ifnum\c@minutes<10 0\fi\the\c@minutes}
\makeatother

\usepackage[OT1]{fontenc}

\usepackage{fancyhdr}
\usepackage{relsize}
\usepackage{bm} 
\usepackage{amssymb}
\usepackage{amsmath}
\usepackage{amsthm}
\usepackage{enumerate}
\usepackage{colortbl}

\usepackage{verbatim} 

\newcommand\skl[1]{{\,\sf#1\,}}
\newcommand\labelcell[1]{\cellcolor[gray]{0.8}\makebox[.5em][c]{\skl{#1}}}
\newcommand\hlinestrut{\hline\rule{0pt}{2.5ex}}
\newcommand\st{$\ast$}
\newcommand\bfgr[1]{\textcolor[gray]{0.5}{\textbf{#1}}}

\newcommand\upla{}
\newcommand\uplapar{}

\newcommand\bnou{\textup{\textbf{\lower3.7pt\hbox{\char'052}:~}}}
\newcommand\enou{\unskip\textup{\textbf{~:\lower3.7pt\hbox{\char'052}}} }
\newcommand\bvell{\textup{\textbf{\lower3.7pt\hbox{\char'052}:~}$\langle$}}
\newcommand\evell{\unskip\textup{$\rangle$\textbf{~:\lower3.7pt\hbox{\char'052}}} }
\newcommand\bbnou{\textup{\textbf{\lower3.7pt\hbox{\char'052\char'052}:~}}}
\newcommand\eenou{\unskip\textup{\textbf{~:\lower3.7pt\hbox{\char'052\char'052}}} }

\newcommand\ie{i.\,e.~}
\newcommand\ifoi{\,\hbox{if\kern2.5pt and\kern2.5pt only\kern2.5pt if}\,{} }

\newcommand\df{\bfseries}
\newcommand\dfc[1]{\,{\df#1}\,}
\newcommand\dfd[1]{\,{\df#1}\hskip1pt}
\newcommand\secpar[1]{\S\,{#1}}

\newcommand\ensep{\unskip\hskip.65em\ignorespaces}

\newcommand\atilde{\smash{\lower.75ex\hbox{\~{}}}}
\newcommand\underl{\smash{\lower.75ex\hbox{-}}}

\newcommand\remark{\medskip\noindent\textit{Remark}.\hskip.5em}

\newcommand\better{$\succ$}
\newcommand\tied{$\sim$}

\newcommand\pq[2]{\raise.25ex\hbox{\footnotesize${#1}\over{#2}$}%
\hskip-.35ex\null}
\newcommand\onehalf{\frac12}

\newcommand\halfsmallskip{\vskip0.5\smallskipamount}

\newcommand\xxxx[1]{%
 \hangindent2.5\parindent
 \hangafter1
 \noindent\hskip.5\parindent
 \hbox to2\parindent{\hss#1\hss}}
\newcommand\condition[2]{\xxxx{#1}\textit{#2}.}
\newcommand\iim[1]{\xxxx{\textup{(#1)}}\ignorespaces}

\newcommand\ddd[1]{\halfsmallskip\vskip-2pt\noindent\hbox to 2\parindent{\hss\footnotesize$\bullet$\ \ }{#1}\ensep}

\newcommand\brwrap[1]{[\textsl{#1}\kern1pt]}

\makeatletter
\newcommand\bibref[1]{\@nameuse{b@#1}}
\renewcommand\@biblabel[1]{\brwrap{#1}}
\renewcommand\@cite[2]{\hbox{\brwrap{#1\if@tempswa\/\upshape\,:\,{\relscale{0.95}#2}\fi}}} 
\makeatother

\newcommand*\dbibref[2]{\bibref{#1}\,\textup{:\,{\relscale{0.95}#2}}}

\newcommand*\refco{\/\kern.1ex\textup{,}\hskip.45ex}
\newcommand*\refsc{\/\kern.15ex\textup{;} }

\newcommand\latop[2]{{#1\atop#2}}

\newcommand\sbset{\subset}

\newcommand\sbseteq{\subseteq}
\newcommand\spseteq{\supseteq}

\newcommand\cd[1]{\!#1\!}

\newcommand\stv{{}^\ast\kern-.25pt v}

\newcommand\isc{v^\ast}
\newcommand\icr{\kappa} 
\newcommand\mast{\kappa} 
\newcommand\img{m^\mast}

\newcommand\rxi{\mathrel{\smash{\succ\kern-1.7ex\raise1.15ex\hbox{\mathsurround0pt$\scriptscriptstyle\xi$}\kern.4ex}}}
\newcommand\rxieq{\mathrel{\smash{%
 \vbox{\offinterlineskip\halign{\hfil##\hfil\cr
 \mathsurround0pt$\succ$\cr
 \noalign{\vskip-.5ex}%
 \mathsurround0pt$-$\cr
 \noalign{\vskip-1.15ex}%
 }\vss}\kern-1.05ex\raise1.15ex\hbox{\mathsurround0pt$\scriptscriptstyle\xi$}\kern.4ex}}}

\newcommand\ppmg{m^\sigma}
\newcommand\ppto{t^\sigma}

\newcommand\psc{v^\pi}
\newcommand\pmg{m^\pi}
\newcommand\pto{t^\pi}
\newcommand\apm{d^\pi} 

\newcommand\vk{v^k}

\newcommand\plumpf{f}

\makeatletter
\gdef\centre#1{\smash{\vbox{\m@th\ialign{\hfil##\hfil\crcr
  \hskip2.25\p@{\tiny$\scriptscriptstyle\bullet$}\crcr
  \noalign{\kern1\p@\nointerlineskip}
  $\hfil\displaystyle{#1}\hfil$\crcr}}}}
\makeatother
\newcommand\pcentre[1]{\left(#1\right)^{\hskip-1.5pt\raise.5pt\hbox{\tiny$\scriptscriptstyle\bullet$}}}

\newcommand\ist{A}
\newcommand\xst{X}
\newcommand\yst{Y}
\newcommand\zst{Z}

\newcommand\tie{\hbox{\textit{\char"05}}}

\newcommand\cst{C}

\newcommand\clustit{\kern.1ex\widetilde c\kern.2ex}

\newcommand\rank[1]{r_{#1}}

\newcommand\rlr[1]{R_{#1}}

\newcommand\vbis{\widetilde v}

\newcommand\iscbis{\widetilde v{}^{\kern.5pt\ast}}
\newcommand\icrbis{\widetilde\icr}
\newcommand\rlrbis[1]{\hbox to1.97ex{\hss\hskip2.5pt$\smash{\widetilde{%
 \hbox to1.9ex{\hss\rule{0pt}{1.45ex}\smash{$R$}\hskip2.5pt\hss}}}$\hss}_{#1}}
\newcommand\pscbis{\widetilde v{}^{\kern.75pt\pi}}

\newcommand\rlrating{\omega}

\newcommand\xibis{\smash{\widetilde{\hbox{\rule{0pt}{1.5ex}\smash{$\xi$}}}}}

\newcommand\ppmgbis{\widetilde m^\sigma}
\newcommand\pptobis{\widetilde t^\sigma}
\newcommand\ptobis{\widetilde t^\pi}

\newcommand\gammabis{\widetilde\gamma}
\newcommand\tref[1]{\setbox0\hbox{\ref{#1}}%
 \smash{\hbox to\wd0{\hss$\widetilde{\rule{0pt}{1.42ex}%
 \smash{\hbox to1em{\hss\ref{#1}\hss}}}$\hss}}}

\newcommand\rev{\hbox{'}}
\newcommand\gpxy{\gamma\hbox{'}\kern-3pt{}_{xy}}
\newcommand\gbispxy{\widetilde\gamma\hbox{'}\kern-3pt{}_{xy}}

\newcommand\gto{\tau}

\newcommand\intg{\gamma}
\newcommand\inth{\eta}

\newcommand\vaa{\mathsf{V}}
\newcommand\vaasub{\vaa\kern-1pt}
\newcommand\vxx{\vaasub_{\scriptscriptstyle X\kern-1pt X}}

\newcommand\vxy{\vaasub_{\scriptscriptstyle X\kern-.25pt Y}}
\newcommand\vyx{\vaasub_{\scriptscriptstyle Y\kern-1pt X}}
\newcommand\vyy{\vaasub_{\scriptscriptstyle Y\kern-.25pt Y}}
\newcommand\vrs{\vaasub_{\scriptscriptstyle R\kern-.25pt S}}
\newcommand\vzz{\vaasub_{\scriptscriptstyle Z\kern-1pt Z}}
\newcommand\vxz{\vaasub_{\scriptscriptstyle X\kern-1pt Z}}

\newcommand\fxx{F_{\scriptscriptstyle X\kern-1pt X}}
\newcommand\fyy{F_{\scriptscriptstyle Y\kern-.25pt Y}}
\newcommand\fxy{F_{\scriptscriptstyle X\kern-.25pt Y}}
\newcommand\fzz{F_{\scriptscriptstyle Z\kern-1pt Z}}
\newcommand\frr{F_{\scriptscriptstyle R\kern-.25pt R}}
\newcommand\flrbis{\hbox to1.75ex{\hss\hskip1.25pt$\widetilde{\hbox to1.5ex{\hss$\varphi$\hskip1.25pt\hss}}$\hss}}
\newcommand\vaabis{\smash{\widetilde{\hbox{\vphantom{t}\smash{$\vaa$}}}}}
\newcommand\vaabissub{\vaabis\kern-1pt}
\newcommand\vaabisxx{\vaabissub_{\scriptscriptstyle X\kern-1pt X}}
\newcommand\vaabiszz{\vaabissub_{\scriptscriptstyle Z\kern-1pt Z}}

\newcommand\xstbiss{\smash{\widetilde\xst}}
\newcommand\ystbiss{\smash{\widetilde\yst}}
\newcommand\zstbiss{\smash{\widetilde\zst}}
\newcommand\vxxbis{\vaabissub_{\scriptscriptstyle \xstbiss\kern-1pt\xstbiss}}
\newcommand\vyxbis{\vaabissub_{\scriptscriptstyle \ystbiss\kern-1pt\xstbiss}}
\newcommand\vzzbis{\vaabissub_{\scriptscriptstyle \zstbiss\kern-1pt\zstbiss}}
\newcommand\flrbissub{\flrbis\kern-1pt}

\newcommand\fbis{\hbox to1.97ex{\hss\hskip2.5pt$\smash{\widetilde{%
 \hbox to1.9ex{\hss\vphantom{t}\smash{$F$}\hskip2.5pt\hss}}}$\hss}}
\newcommand\fbisxx{\fbis_{\kern-2pt\scriptscriptstyle \xstbiss\kern-1pt\xstbiss}}

\newcommand\av[1]{\alpha_{#1}}

\newcommand\vbisbis{\smash{\widetilde{\hbox{\vphantom{t}\smash{$\vbis$}}}}\kern.5pt}
\newcommand\iscbisbis{\smash{\widetilde{\hbox{\vphantom{t}\smash{$\vbis$}}\kern.75pt}}^\ast}
\newcommand\icrbisbis{\smash{\widetilde{\hbox{\vphantom{t}\smash{$\icrbis$}}}}\kern.5pt}

\newcommand\chrela{\mathrel{\hbox{\relscale{1.25}$\kern1pt\triangleright\kern1pt$}}}

\newcommand\chrelabis{\mathrel{\smash{\widetilde{\hbox{\vrule width0pt height6.5pt\smash{\hbox{\relscale{1.25}$\kern1pt\triangleright\kern1pt$}}}}}}}

\newcommand\funbis{{\cal F}^{\kern.5pt\prime}}

\newcommand\pptoset{{\mathcal T}\hskip-1pt}
\newcommand\pptosetbis{\smash{\widetilde{\mathcal T}}\hskip-1pt}

\newtheorem{proposition}{Proposition}[section]
\newtheorem{lemma}[proposition]{Lemma}
\newtheorem{theorem}[proposition]{Theorem}
\newtheorem{corollary}[proposition]{Corollary}

\newlength\repskip 
{\null\hskip-\repskip\hbox to0.9\hsize\bgroup\hbox to0pt{\small$(\ref{#1})$\hss}\hfil\hskip.1\hsize$\displaystyle}%
{$\hfil\egroup}

{\hbox to\hsize\bgroup\hbox to0pt{\small$(\ref{#1})$\hss}\hfil$\displaystyle}%
{$\hfil\egroup\nonumber}

\begin{document}


\thispagestyle{empty}

\null\vskip-5.5\baselineskip\null 
\renewcommand\upla{\vskip-6pt}

\begin{center}
\hrule
\upla
\vskip7.5mm
\textbf{\uppercase{A continuous rating method for~preferential~voting.\ \ The incomplete case}}
\par\medskip
\textsc{Rosa Camps,\, Xavier Mora \textup{and} Laia Saumell}
\par
Departament de Matem\`{a}tiques,
Universitat Aut\`onoma de Barcelona,
Catalonia,
Spain
\par\medskip
\texttt{xmora\,@\,mat.uab.cat}
\par\medskip
July 13, 2009;\ensep revised October 10, 2011
\vskip7.5mm
\upla
\hrule
\end{center}

\upla
\vskip-4mm\null
\begin{abstract}
A method is given for quantitatively rating the social acceptance of different options which are the matter of a preferential vote.\linebreak
In~contrast to a previous article, here the individual votes are allowed to be incomplete,
that is, they need not express a comparison between every pair of options.
This includes the case where each voter gives an ordered list  
restricted to a subset of most preferred options.
In this connection, the proposed method
(except for one of the given variants)
carefully distinguishes a lack of information
about a given pair of options from a proper tie between them.
As in the special case of complete individual votes,
the~proposed generalization is proved to have certain desirable properties, which include:\ensep
the conti\-nuity of the rates with respect to the data,\ensep
a~decomposition property that characterizes certain situations opposite to a tie,\ensep
the Condorcet-Smith principle,\ensep 
and
clone consistency.

\upla

\upla
\bigskip\noindent
\textbf{AMS subject classifications:} 
\textit{%
05C20, 
91B12, 
91B14, 
91C15, 
91C20. 
}
\end{abstract}

\vskip5mm
\upla
\hrule

\vskip-12mm\null 
\section*{}
In a previous article \cite{crc} we introduced a method for quantitatively rating the social acceptance of different candidate options which are the matter of a preferential vote.
The quantitative character of this method lies in a combination of two properties:
First, a property of continuity that allows to sense the closeness between two candidates,
for instance the winner and the runner-up. 
Second, a property of decomposition that allows to recognise certain situations that are opposite to a tie;
in particular, a candidate gets the best possible rate \ifoi it is placed first by all voters.
\ensep
These two properties are combined with other ones of a qualitative nature.
Especially outstanding among the latter is 
the Condorcet-Smith principle:
Assume that the set of candidates is partitioned into two classes~$\xst$ and~$\yst$ such that for each member of $\xst$ and every member of~$\yst$ there are more than half of the individual votes where the former is preferred to the latter; in that case, the social ranking prefers also each member of $\xst$ to any member of $\yst$.
\ensep
To our knowledge, the existing literature does not offer any other rating method that combines the three mentioned properties, namely continuity, decomposition and the Condorcet-Smith principle.

However, in \cite{crc} we restricted ourselves to the complete case, that is, we assumed that every individual expresses a~comparison (a~preference or a~tie) about each pair of options.
Such a restriction leaves out many cases of interest,
like for instance 
truncated rankings.
In this article we will extend that method to the incomplete case,
where the individual votes need not express a comparison about every pair of options.
This extension will be done in such a way that 
a lack of information about the preference of a voter between a given pair of options 
will be carefully distinguished from a proper tie between them.
\ensep
In this connection, we will have to cope with the fact that a quantitative specification of the collective opinion about a pair of options
has then two degrees of freedom; in fact, knowing how many voters preferred $x$ to $y$ does not determine how many of them preferred $y$ to $x$.
This introduces a special difficulty that was not present in the complete case.


An extreme case of incompleteness is that where each vote reduces to choosing a single option. 
In that case our rates will be linearly related to the vote fractions.
\ensep
Another important case is that of approval voting. In~that case, our~rates will be essentially different from the number of received approvals;
however, one of the variants of our method, namely the margin-based variant, will be shown to rank the options exactly in the same way as the number of received approvals.


\renewcommand\upla{\vskip-1.5pt}

\upla
\medskip
We call our quantitative method the \dfd{CLC~rating method}, where the capital letters stand for ``Continuous Llull Condorcet''.
\ensep
The reader interested to try it can use the \,\textsl{CLC~calculator}\, which has been made available at~\cite{clc}.
\ensep
Of course, any rating automatically implies a ranking.
In this connection, it should be noticed that
the CLC~rating method is built upon certain variants of the qualitative ranking method that was introduced in~1997 by~Markus Schulze
\hbox{\brwrap{
\bibref{sc}\refco
\bibref{scbis}\refsc
\dbibref{t6}{p.\,228--232}%
}}.




\upla
\medskip
This article is organized as follows: Section~1 gives a more precise statement of the problem together with some general remarks. Section~2 presents an outline of the proposed method, followed by a discussion of certain special cases, a summary of the procedure, and a discussion of certain variants. Section~3 gives three representative examples. Finally, sections~\ref{sec-projection}--\ref{sec-approval} give detailed mathematical proofs of the claimed properties.

\renewcommand\upla{}

\section{Statement of the problem and general remarks}

\paragraph{1.1}
Let us consider a set of $N$~options which are the matter of a vote.
Let us assume that each voter expresses his preferences about certain pairs of options.
Our aim is to combine such individual preferences so as to rate the social acceptance of each option on a continuous scale.
More specifically, we would like to do it in accordance with
the following conditions:

\smallskip
\newcommand\llsi{\textup{A}}
\condition{\llsi}{Scale invariance (homogeneity\hskip.1pt)}
The rates depend only on the relative frequency of each possible content of an individual vote.
In~other words, if every individual vote is replaced by a fixed number of copies of it, the rates remain exactly the same.

\smallskip
\newcommand\llpe{\textup{B}}
\condition{\llpe}{Permutation equivariance (neutrality\hskip.1pt)}
Applying a certain permutation of the options to all of the individual votes has no other effect than getting the same permutation in the social rating.

\smallskip
\newcommand\llco{\textup{C}}
\condition{\llco}{Continuity}
The rates depend continuously on the relative frequency of each possible content of an individual vote.

\newcommand\condr[1]{#1\rlap{\mathsurround0pt$_r$}}
\newcommand\condf[1]{#1\rlap{\mathsurround0pt$_f$}}

\medskip
The next two conditions consider a specific form of rating.
From now on, we will refer to it as rank-like rating.

\smallskip
\newcommand\llrr{\textup{D}}
\condition{\llrr}{Rank-like form}
Each rank-like rate is a number, integer or fractional, between $1$ and~$N$. The~best possible value is~$1$ and the worst possible one is~$N$. The average rank-like rate is larger than or equal to~$(N+1)/2$, with equality in the complete case.

\smallskip
\newcommand\llrd{\textup{E}}
\condition{\llrd}{Rank-like decomposition in the complete case}
Assume that the individual preferences are complete, \ie a comparison (a preference or~a~tie) is expressed about every pair of options.
Consider a splitting of the options in two classes~$\xst$ and~$\yst$ such that each member of~$\xst$ is unanimously preferred to every member of~$\yst$.
Such a situation translates into the three following equivalent facts, 
where $|\xst|$~denotes the number of elements of $\xst$:
\ensep
(a)~The rank-like rates of~$\xst$ coincide with those that one obtains when the individual votes are restricted to~$\xst$.
\ensep
(b)~After diminishing them by the number~$|\xst|$, the rank-like rates of~$\yst$ coincide with those that one obtains when the individual votes are restricted to~$\yst$.
\ensep
(c)~The average rank-like rate of~$\xst$ is $(|\xst|+1)/2$.



\smallskip
\noindent
In particular, in the complete case an option will get a rank-like rate exactly equal to~$1$ [resp.~$N$]
\,if and only if\, it is unanimously preferred to [resp.~considered worse than] any other.
As we will see later on, some of the implications contained in the preceding condition
will hold also in certain situations that allow for incompleteness.


\medskip
The next condition considers the extreme case of incompleteness where each vote reduces to choosing a single option.

\smallskip
\newcommand\llrp{\textup{F}}
\condition{\llrp}{Rank-like rates for single-choice voting}
Assume that each vote reduces to choosing a single option. In that case, the rank-like rate of each option is the weighted average of the numbers $1$ and $N$ with weights given respectively by the fraction of the vote in favour of that option \,and\, the complementary fraction.

\medskip
Finally, we require a condition that concerns only
the concomitant social ranking, that is,
the ordinal information contained in the social rating:

\smallskip
\newcommand\llmp{\textup{M}}
\condition{\llmp}{Condorcet-Smith principle}
Consider a splitting of the options into two classes~$\xst$ and~$\yst$.
Assume that for each member of $X$ and every member of $Y$ there are more than half of the individual votes where the former is preferred to the latter.
In that case, the social ranking also prefers each member of $X$ to every member of~$Y$.


\paragraph{1.2}
We will adopt the point of view of paired comparisons.
In other words, we will be based upon the number of voters who prefer $x$ to $y$,
where $x$ and $y$ vary over all ordered pairs of options.
These numbers will be denoted by $V_{xy}$. 
Most of the time, however, we will be dealing with the fractions $v_{xy}=V_{xy}/V$,
where $V$ denotes the total number of votes.
We will refer to $V_{xy}$ and $v_{xy}$ respectively
as the absolute and relative \dfc{scores} of the pair $xy$,
and the whole collection of these scores 
will be called the (absolute or relative)
\dfc{Llull matrix} of the vote.

In the complete case one has 
$v_{xy}+v_{yx}=1$, 
so $v_{xy}$ automatically determines $v_{yx}$.
In contrast, in the incomplete case we are ensured only that $v_{xy}+v_{yx}\le1$,
so $v_{xy}$ alone does not determine $v_{yx}$.
\ensep
In particular, the conditions \,$v_{xy}>\onehalf$ \,and\, $v_{xy}>v_{yx}$\, are not equivalent to each other, which gives rise to two possible notions of majority.
\ensep
Anyway, in the incomplete case a quantitative specification of preference between two options $x$~and~$y$ requires the values of both $v_{xy}$ and $v_{yx}$,
or equivalently, 
their sum and difference, 
$t_{xy} = v_{xy}+v_{yx}$ and $m_{xy} = v_{xy}-v_{yx}$, which we will call respectively the (relative) \dfc{turnout} and \dfc{margin} associated with the pair $xy$.

\paragraph{1.3}
In preferential voting
the individual preferences are usually assumed to be expressed in the form of a ranking,
that is, a list of options in order of preference.
In this connection, it is quite natural to admit the possibility of ties as well as incomplete lists.
When we are dealing with incomplete lists, their translation into paired comparisons admits of several interpretations.
\ensep
In most cases, it is reasonable to use the following one:

\iim{a}When $x$ and $y$ are both in the list\,
and $x$ is ranked above $y$ (without a tie),
we certainly interpret that $x$ is preferred to $y$.

\iim{b}When $x$ and $y$ are both in the list\,
and $x$ is ranked as good as $y$,\,
we interpret it as being equivalent to
half a vote preferring $x$ to $y$\,
plus another half a vote preferring~$y$~to~$x$.

\iim{c}When $x$ is in the list and $y$ is not in it,\,
we interpret that $x$ is preferred to $y$.

\iim{d}When neither $x$ nor $y$ are in the list,\,
we interpret nothing
about the preference of the voter between $x$ and $y$.

\medskip
Instead of rule~(d), one can consider the possibility of using the following alternative:

\iim{d$'$} When neither $x$ nor $y$ are in the list,\,
we interpret that they are considered equally good
(or equally bad),\, so we proceed as in~(b).

\noindent
This amounts to complete each truncated ranking by appending to it all the missing options tied to each other, which brings the problem to the complete case considered in~\cite{crc}.
\ensep
Generally speaking, however, this interpretation can be criticized in that the added information might not be really meant by the voter.

\medskip
On the other hand, in the spirit of not adding any information not really meant by the voter, in certain cases it may be appropriate to replace rule~(c) by the following one:

\iim{c$'$} When $x$ is in the list and $y$ is not in it,\,
we interpret nothing about the preference of the voter between $x$ and $y$.

\medskip
Generally speaking, 
the individual votes could be arbitrary binary relations,
interpreted as it is mentioned in \cite[\secpar{3.1}]{crc};
even more generally, they could be valued binary relations belonging to~$\Omega\!=\!\{v\in[0,1]^{\textit{\char"05}} \mid v_{xy}+v_{yx}\le1\}$, where $\tie$~denotes the set of pairs $xy\in\ist\times\ist$ with $x\neq y$ \cite[\secpar{3.3}]{crc}.
Such a~possibility makes sense in that the individual opinions may already be the result of aggregating a variety of criteria.

Anyway, the collective Llull matrix is simply the center of gravity of a~distribution of individual votes:
\begin{equation}
\label{eq:cog}
v_{xy} = \sum_k \alpha_k\, \vk_{xy},
\end{equation}
where $\alpha_k$ are the relative frequencies or weights of the individual votes~$\vk\!\in\!\Omega$.

\paragraph{1.4}
In the particular case where the set $\xst$ consists of a single option, the~Condorcet-Smith principle~\llmp\,\ takes the following form:

\smallskip
\newcommand\llmpw{\textup{\llmp 1}}
\condition{\llmpw}{Condorcet principle (majority form)} If an option $x$ has the property that for every $y\neq x$ there are more than half of the individual votes where $x$ is preferred to $y$, then $x$ is the social winner.

\smallskip
\noindent
In the complete case (where the Condorcet principle was originally proposed) the preceding condition is equivalent to the following one:

\smallskip
\newcommand\llcpw{\textup{\llmp 1$'$}}
\condition{\llcpw}{Condorcet principle (margin form)} If an option $x$ has the property that for~every $y\neq x$ there are more individual votes where $x$ is preferred to $y$\, than vice versa,\, then $x$ is the social winner.

\smallskip
\noindent
However,
generally speaking condition~\llmpw\ is weaker than~\llcpw,
and the CLC~me\-thod will satisfy only the weaker version.

\smallskip
This~lack of~compliance with the stronger condition \llcpw\ 
may be considered undesirable.
\ensep
However, other authors have already remarked the need to 
require only \llmpw\ 
in order to be able to keep other properties (see for instance~\cite{wo}).
In our case, \llcpw\ seems to conflict with the continuity property~\llco.

\paragraph{1.5}
Conditions \llrd\ and \llmp\ refer to cases where all the voters or at least half of them proceed in a certain way.
Of course, it should be clear whether all the voters means all of the actual ones or maybe all the potential ones (\ie actual voters plus abstainers).
We assume that one has made a choice in that connection,
thus defining the total number of voters~$V$.
Mathematically speaking, we only need $V$ to be larger than any  absolute turnout~$V_{xy}+V_{yx}$.
Increasing the value of $V$ 
has no other effect than contracting
the final rating towards the point where all rates take the worst possible value ($N$~for rank-like rates).


\section{Outline of the method}

This section presents the proposed method at the same time that it introduces the associated terminology.
As in \cite{crc}, the procedure involves 
a projection of the Llull matrix onto a special set of such matrices.
\ensep
Steps~0--3, as well as the last one, 
will be exactly the same as in the complete case.
However, steps~4 and 5 contain new elements.
More specifically, step~4 
requires a quadratic optimization in connection with the turnouts,
and step~5 takes the union of certain intervals
where the complete case takes the maximum of certain margins.
\ensep
The reasons behind steps 1--3 and 6 were explained in \cite[\secpar{2}]{crc}.
Those behind steps~4 and 5 will be briefly explained in \secpar{2.2.3}, after having looked at the particularities of the complete case as well as those of single-choice voting.


\newcommand\step[1]{\bigskip\noindent\textit{Step~#1.}\hskip.5em}

\paragraph{2.1}
\textit{Step~0.}\hskip.5em
To begin with, we must form the Llull matrix~$(v_{xy})$.
Its entries are the relative \dfc{scores} $v_{xy} = V_{xy}/V$,
where $V$ is the number of voters, 
and $V_{xy}$ counts how many of them prefer~$x$ to $y$.
In the case of 
ranking votes, 
this count will usually make use of rules~(a--d) of \secpar{1.3},
though in certain cases it may be reasonable to replace rule (c) by (c$'$), or rule (d) by (d$'$).
Besides the scores themselves, which will be used in the next step,
later on we will also make use of the associated \dfc{turnouts}
\begin{equation}
\label{eq:turnout}
t_{xy} = v_{xy} + v_{yx}.
\end{equation}



\step{1} 
Concerning the margin component,
we will rely on the \dfc{indirect scores}~$\isc_{xy}$. They derive from the original scores through an operation that generalizes the notion of transitive closure to valued relations.
More specificaly, they are defined in the following way:
\begin{equation}
\label{eq:isc}
\isc_{xy} \,=\, \max\,\{v_\alpha\mid \text{$\alpha$ is a path \,$x_0
x_1 \dots x_n$\, from $x_0=x$ to $x_n=y$}\,\},
\end{equation}
where the score~$v_\alpha$ of a path $\alpha=x_0 x_1 \dots x_n$ is
defined as
\begin{equation}
\label{eq:valpha}
v_\alpha \,=\, \min\,\{v_{x_ix_{i+1}}\mid 0 \le i < n\,\}.
\end{equation}
Obviously, we have 
\begin{align}
0 \,&\le\,\isc_{xy}\le 1,
\\[2.5pt]
&\isc_{xy} \,\ge\, v_{xy}.
\end{align}
The indirect turnouts $\isc_{xy}+\isc_{yx}$ can be larger than $1$.
However, the following steps will use the indirect scores only through the associated \dfc{indirect margins}
\begin{equation}
\label{eq:img}
\img_{xy}\,=\,\isc_{xy}-\isc_{yx}.
\end{equation}

\smallskip
\step{2} 
This step is the discrete core of the procedure.
It begins by considering the \dfc{indirect comparison relation}
\begin{equation}
\label{eq:nu}
\icr \,=\, \{\,xy\mid  \img_{xy} > 0\,\},
\end{equation}
as well as its codual
\begin{equation}
\label{eq:hatnu}
\hat\icr \,=\, \{\,xy\mid  \img_{xy} \ge 0\,\}.
\end{equation}
The relation~$\icr$ has the virtue of being transitive.
This crucial fact was remarked in 1998 by Markus Schulze~%
\hbox{\brwrap{\bibref{sc}\,b}}. 
As a consequence, $\icr$ is a partial order, and therefore 
one can always extend it to a total order~$\xi$.
\ensep
For instance, according to Proposition~5.2 of \cite{crc},
it suffices to arrange the options by non-decreasing values of the ``tie-splitting'' Copeland ranks
\begin{equation}
\label{eq:copeland}
\rank x \,\,=\,\, 1 \,+\, |\{\,y\mid y\cd\neq x,\ \img_{yx}\cd>0\}|
\,+\, {\hbox{\large$\frac12$}}\,|\{\,y\mid y\cd\neq x,\ \img_{yx}\cd=0\}|.
\end{equation}
Such a total order extension $\xi$ automatically satisfies
not only $\icr\sbseteq\xi$ but also $\xi\sbseteq\hat\icr$.
In the following we call it an \dfc{admissible order}.

The following steps 
assume that one has fixed an 
admissible order~$\xi$.
\ensep
The intermediate quantities computed in these steps may depend on which admissible order is used, but the final results will be independent of it.
From now on,
the situation $xy\in\xi$ will be expressed also by writing $x\rxi y$. 
Furthermore, $x'$ will mean the immediate successor of $x$ in~$\xi$.

\step{3} 
Starting from the indirect margins $\img_{xy}$, one computes the \dfc{super\-diagonal inter\-mediate projected margins}
\begin{equation}
\label{eq:ppmg}
\ppmg_{xx'} \,=\, \min\,\{\, \img_{pq} \;\vert\; p\rxieq x,\; x'\rxieq q\,\}.
\end{equation}


\step{4} 
Starting from the original turnouts $t_{xy}$, one computes the \dfc{inter\-mediate projected turnouts}~$\ppto_{xy}$, These numbers depend not only on the original turnouts $t_{xy}$, but also on the superdiagonal intermediate projected margins $\ppmg_{xx'}$. More specifically, they are 
taken as the values of $\gto_{xy}$ that minimize the quantity
\begin{equation}
\Phi \,=\, \sum_x \, \sum_{y\neq x} \, (\gto_{xy} - t_{xy})^2
\label{eq:phi}
\end{equation}
under the following constraints:
\begin{gather}
\gto_{xy} \,=\, \gto_{yx},
\label{eq:symmetry}
\\[2.5pt]
\ppmg_{xx'} \,\le\, \gto_{xx'}\,\le\, 1,
\label{eq:bounded}
\\[2.5pt]
0 \,\le\, \gto_{xz} - \gto_{x'z} \,\le\, \ppmg_{xx'},
\quad\hbox{whenever $z\not\in\{x,x'\}$.}
\label{eq:increment}
\end{gather}
The actual computation of the minimizer can be carried out in a finite number of steps by means of a quadratic programming algorithm
\cite[ch.\,16]{nw}.
\ensep
For future reference, the set of matrices $(\gto_{xy})$ that satisfy (\ref{eq:symmetry}--\ref{eq:increment}) will be denoted as~$\pptoset$, \,and the preceding minimizing operation that defines the intermediate projected 
turnouts~$(\ppto_{xy})$ as a function of the original turnouts~$(t_{pq})$ and the superdiagonal 
intermediate projected margins~$(\ppmg_{pp'})$ will be denoted as~$\Psi$:
\begin{equation}
\label{eq:psi}
\ppto_{xy} \,=\, \Psi[(t_{pq}),(\ppmg_{pp'})]_{xy}.
\end{equation}

\step{5} 
Form the intervals
\begin{equation}
\label{eq:gammaxxp}
\gamma_{xx'} \,=\, [\,(\ppto_{xx'}-\ppmg_{xx'})/2\,,\, (\ppto_{xx'}+\ppmg_{xx'})/2\,],
\end{equation}
as well as their unions
\begin{equation}
\label{eq:gammaxy}
\gamma_{xy} \,=\, \bigcup\,\,\gamma_{pp'}, \,\,\,\hbox{with \,$p$\, varying in the interval \,$x\rxieq p\rxi y$,\hfil}
\end{equation}
where $xy$ is restricted to satisfy $x\rxi y$. The sets $\gamma_{xy}$ are still intervals.
The~\dfc{projected scores} are the upper and lower bounds of these intervals:
\begin{equation}
\label{eq:psc}
\psc_{xy} \,=\, \max\,\gamma_{xy},\qquad
\psc_{yx} \,=\, \min\,\gamma_{xy}.
\end{equation}

\smallskip
Equivalently to (\ref{eq:gammaxxp}--\ref{eq:psc}),
the projected scores can be computed also in the following way,
which has a more practical character:
Take
\begin{equation}
\label{eq:pscxxp}
\psc_{xx'} = (\ppto_{xx'}+\ppmg_{xx'})/2,\qquad
\psc_{x'x} = (\ppto_{xx'}-\ppmg_{xx'})/2,
\end{equation}
and then, for $x\rxi y$,
\begin{equation}
\label{eq:pscxy}
\psc_{xy} = \max\,\{\psc_{pp'} \,\vert\, x\rxieq p\rxi y\},\qquad
\psc_{yx} = \min\,\{\psc_{p'p} \,\vert\, x\rxieq p\rxi y\}
\end{equation}

\step{6} 
Finally, the rank-like rates are determined by the formula
\begin{equation}
\label{eq:rrates}
\rlr{x} \,=\, N - \sum_{y\ne x} \psc_{xy}
\end{equation}

\renewcommand\uplapar{\vskip-8mm\null}

\paragraph{2.2}
\textbf{Special cases and heuristic considerations.}\hskip.5em

\uplapar
\subparagraph{2.2.1}
In the complete case the original turnouts $t_{xy}$ are all of them equal to~$1$.
One easily sees that in these circumstances step~4 results in $\ppto_{xy}$ also equal to~$1$ for all pairs~$xy$. In fact, this choice clearly satisfies conditions~(\ref{eq:symmetry}--\ref{eq:increment}) at the same time that it certainly minimizes (\ref{eq:phi}) to $0$.
As a consequence, the intervals $\gamma_{xx'}$ defined by (\ref{eq:gammaxxp}) are all of them centred at~$1/2$. Obviously, this property will be inherited by their unions $\gamma_{xy}$. On the other hand, one easily sees that the width of $\gamma_{xy}$, in other words the projected margin $\pmg_{xy}=\psc_{xy}-\psc_{yx}$, will be the following: $\pmg_{xy} = \max\,\{\ppmg_{pp'} \,\vert\, x\rxieq p\rxi y\}$.\linebreak 
So,~\textit{in the complete case the above-described procedure reduces to the one that was presented in}~\cite{crc}.

\uplapar
\subparagraph{2.2.2}
Let us see what we get in the case of single-choice voting.
\ensep
To~begin with, rules~(c--d) of \secpar{1.3} result in $v_{xy} = \plumpf_x$ for~every $y\neq x$, where $\plumpf_x$ is the fraction of voters who choose~$x$. 
\ensep
This implies that $\isc_{xy} = v_{xy} = \plumpf_x$. In fact, any path $\gamma$ from $x$ to $y$ starts with a link of the form $xp$, whose associated score is $v_{xp} = \plumpf_x$. So $v_\gamma \le \plumpf_x$ and therefore $\isc_{xy} \le \plumpf_x$. But on the other hand $\plumpf_x = v_{xy} \le \isc_{xy}$.
\ensep
Consequently, we get $\img_{xy} = \isc_{xy} - \isc_{yx} = v_{xy} - v_{yx} = \plumpf_x - \plumpf_y$, and the admissible orders are those for which the $\plumpf_x$ are non-increasing.
\ensep
Owing to this non-increasing character, the intermediate projected margins are~$\ppmg_{xx'} = m_{xx'} = \plumpf_x - \plumpf_{x'}$.
\ensep
On the other hand, the intermediate projected turnouts are~$\ppto_{xy} = t_{xy} = \plumpf_x + \plumpf_{y}$. In fact these numbers are easily seen to satisfy conditions (\ref{eq:bounded}--\ref{eq:increment}) and they obviously minimize~(\ref{eq:phi}).
\ensep
As a consequence, $\gamma_{xx'} = [\plumpf_{x'},\plumpf_x]$. In particular, the intervals $\gamma_{xx'}$ and $\gamma_{x'x''}$ are adjacent to each other (the right end of the latter coincides with the left end of the former). This fact entails that $\gamma_{xy} = [\plumpf_y,\plumpf_x]$ whenever $x\rxi y$.
\ensep
So,~the projected scores are the end points of these intervals, namely $\psc_{xy} = \plumpf_x$ and $\psc_{yx} = \plumpf_y$.
In particular, they coincide with the original scores.
\ensep
Finally, \textit{the rank-like rates are $\rlr x=1+(N-1)f_x=f_x + (1-f_x)N$, as stated by condition~\llrp}.



\uplapar
\subparagraph{2.2.3}
In \secpar{2.2.1} we have seen that in the complete case, the projected margins are obtained from the superdiagonal intermediate ones by means of a maximum operation: $\pmg_{xy} = \max\,\{\ppmg_{pp'} \,\vert\, x\rxieq p\rxi y\}$ (whenever $x\rxi y$).\linebreak 
In contrast, in the case of single-choice voting we have a sum: $\pmg_{xy} = \sum\,\{\ppmg_{pp'} \,\vert\, x\rxieq p\rxi y\}$ (whenever $x\rxi y$, since we have both $\pmg_{xy} = m_{xy} = \plumpf_x-\plumpf_y$ and $\ppmg_{pp'} = m_{pp'} = \plumpf_p-\plumpf_{p'}$); as we have just seen in \secpar{2.2.2}, it must necessarily be so if we are to satisfy condition~\llrp. So a general method requires an operation that reduces to maximum in one case and to addition in the other.

This leads to the idea that this general operation should be that of taking the union of suitable ``score intervals''. A~score interval can be viewed as giving a pair of scores about two options, these scores being respectively in favour and against a specified preference about the two options. Alternatively, it can be viewed as giving a certain margin together with a certain turnout.
Following on this line, one can be tempted to also replace the minimum operation (\ref{eq:ppmg}) of step~3 by an intersection of the score intervals that combine the original turnouts~$t_{xy}$ with the indirect margins~$\img_{xy}$. 
Such a procedure works as desired both in the case of complete votes and that of single-choice voting. However, it breaks down in other cases of incomplete votes that produce empty intersections and disjoint unions.

In order to avoid these problems, we were led to replace the original turnouts $t_{xy}$ by the intermediate projected ones $\ppto_{xy}$.
As we will see, the constraints that are imposed on the latter have the virtue of ensuring a non-empty intersection for any two consecutive intervals $\gamma_{xx'}$ and $\gamma_{x'x''}$. 
On the other hand, they ensure also the inequality $\pto_{xx'}\ge\pto_{x'x''}$.
These facts will be crucial for achieving 
the following properties for the final rank-like rates $\rlr x$:\ensep
(a)~being the same for any admissible order~$\xi$;\ensep
and (b)~being consistent with any such order $\xi$, \ie having $\rlr x \le \rlr y$ whenever $x\rxi y$.


\renewcommand\uplapar{\vskip-6.5mm\null}

\uplapar
\paragraph{2.3}
\textbf{Summary}
\uplapar
\begin{enumerate}
\setlength\itemsep{0pt}
\setcounter{enumi}{-1}
\item Form the Llull matrix~$(v_{xy})$. Work out the turnouts~$t_{xy}=v_{xy}+v_{yx}$.
\item Compute the indirect scores~$\isc_{xy}$ defined by~(\ref{eq:isc}--\ref{eq:valpha}). An efficient way to do it is the Floyd-Warshall algorithm \cite[\secpar{25.2}]{co}.
For small values of $N$, one can do it by hand in successive steps that progressively increase the length of the paths under consideration.
Having computed the indirect scores, one works out 
the indirect margins $\img_{xy} \cd= \isc_{xy}\cd-\isc_{yx}$.
\item Consider the indirect comparison relation $\icr = \{xy\mid\img_{xy} > 0\}$. 
Fix an admissible order~$\xi$, \ie a total order that extends $\icr$.
For instance, it suffices to arrange the options by non-decreasing values of the ``tie-splitting'' Copeland ranks (\ref{eq:copeland}).
\item Starting from the indirect margins~$\img_{xy}$, work~out the superdiagonal intermediate projected margins $\ppmg_{xx'}$ as defined in (\ref{eq:ppmg}).
\item Starting from the original turnouts~$t_{xy}$, and taking into account the superdiagonal intermediate projected margins $\ppmg_{xx'}$, determine the intermediate projected turnouts $\ppto_{xy}$ by minimizing (\ref{eq:phi}) under the\linebreak constraints (\ref{eq:symmetry}--\ref{eq:increment}). This can be carried out in a finite number of steps by means of a quadratic programming algorithm
\cite[ch.\,16]{nw}.
\item Form the intervals $\gamma_{xx'}$ defined by (\ref{eq:gammaxxp}), derive their unions $\gamma_{xy}$ as defi\-ned by~(\ref{eq:gammaxy}),\, and read off the projected scores~$\psc_{xy}$ (\ref{eq:psc}).
\vskip-3.5ex\null
\item[] \hfil Or equivalently:
\vskip-3.5ex\null
\item[] Compute the sub- and super-diagonal projected scores
as defined by~(\ref{eq:pscxxp}), and then derive all the others according to~(\ref{eq:pscxy}).

\item Compute the rank-like rates~$\rlr{x}$ according to (\ref{eq:rrates}).
\end{enumerate}

\vskip-3mm
\uplapar
\paragraph{2.4}
\textbf{Variants.}\hskip.5em
The preceding procedure admits of certain variants which might be appropriate to some special situations. Next we will distinguish four of them, namely

\newcommand\iitm[2]{%
\halfsmallskip\hskip2\parindent\hbox to1.5em{#1.\hss}#2}

\iitm{1}{Main}
\iitm{2}{Codual}
\iitm{3}{Balanced}
\iitm{4}{Margin-based}
\smallskip

\renewcommand\upla{}
\renewcommand\uplapar{}

\noindent
The above-described procedure is included in this list as the main variant. The four variants are exactly equivalent to each other in the complete case, but in the incomplete case they can produce different results. In spite of this, they all share the main properties.
Having said that, the proofs given in this paper assume either the main variant or the margin-based one.

\smallskip
The codual variant is analogous to the main one except that the \hbox{max-min} indirect scores $\isc_{xy}$ are replaced by the following \hbox{min-max} ones:
\begin{equation}
\stv_{xy} \hskip.75em = \hskip.75em
\min_{\vtop{\scriptsize\halign{\hfil#\hfil\cr\noalign{\vskip.5pt}$x_0=x$\cr$x_n=y$\cr}}}
\hskip.75em
\max_{\vtop{\scriptsize\halign{\hfil#\hfil\cr\noalign{\vskip-1.25pt}$i\ge0$\cr$i<n$\cr}}}
\hskip.75em v_{x_ix_{i+1}}.
\label{eq:codual}
\end{equation}
Equivalently, $\stv_{xy} = 1 - {\hat v}^\ast_{yx}$ where ${\hat v}_{xy} = 1 - v_{yx}$. In the complete case one has $\stv_{xy} = 1 - \isc_{yx}$, so that $\stv_{xy} - \stv_{yx} = \isc_{xy} - \isc_{yx}$;
as a consequence, the codual variant is then equivalent to the main one.

\smallskip
The balanced variant takes
$\icr = \{xy\mid \isc_{xy} \cd> \isc_{yx}, \stv_{xy} \cd> \stv_{yx}\}$
together with
\begin{equation}
\img_{xy} = \begin{cases}
\min\,(\isc_{xy} - \isc_{yx}, \stv_{xy} - \stv_{yx}),
& \text{if $xy\in\icr$,}\\
-\img_{yx},
& \text{if $yx\in\icr$,}\\
0,
& \text{otherwise.}\\
\end{cases}
\label{eq:balanced}
\end{equation}
The remarks made in connection with the codual variant show that in the complete case the balanced variant is also equivalent to the preceding ones.

\smallskip
The margin-based variant follows the procedure of \cite{crc} even if one is not originally in the complete case.
Equivalently, it corresponds to replacing the original scores $v_{xy}$ by the following ones: $v'_{xy} = (1+m_{xy})/2$, where $m_{xy}=v_{xy}-v_{yx}$.
This amounts to replacing any lack of information about a pair of options by a proper tie between them, which brings the problem into the complete case.
In the case of ranking votes, it corresponds to interpreting the votes using rule (d$'$) instead of rule (d) as mentioned in~\secpar{1.3}.
So, the specific character of this variant lies only in its interpretation of incomplete votes.

\remark
Other variants ---in the incomplete case--- arise when equation (\ref{eq:rrates}) is replaced by the following one:
\begin{equation}
\rlr{x} \,=\, 1 + \sum_{y\neq x}\,\psc_{yx}.
\label{eq:rratesVar}
\end{equation}

\paragraph{2.5}
The main ideas that underlie the indirect scores of step~1 and the associated indirect comparison relation $\icr$ are the same as in the 
ranking method 
of~Markus Schulze \cite{scbis}.
Like us, he allows for the Llull matrix to be incomplete,
and he distinguishes several variant methods.
One of them is a~margin variant that coincides essentially with ours;
more specifically, both of them give exactly the same 
indirect comparison relation~$\icr$ ---Schulze's~$\mathcal{O}$--- 
although they may differ in the subsequent treatment of ties.

However, aside from the margin approach, which does not properly face incompleteness,
none of our other variants coincides with any of Schulze's ones.
As a matter of fact, they can result
in different indirect comparison relations~$\icr$ 
and therefore, by Theorem~\ref{st:RvsNu} below, in different final rankings.

In this connection, we may add that none of Schulze's variants except the margin one 
is appropriate for being extended to a continuous rating method.
In fact, to this effect, the strength comparison relation $\succ_D$ that Schulze
uses to compare the pairs $(v_{xy}, v_{yx})$ with each other 
should remain unchanged under small perturbations of the scores.
However, this is not the case for any of those variants,
namely {\small$D=\text{ratio}$}, {\small$D\,=\,\text{win}$},
and {\small$D=\text{los}$}.
More specifically, by looking at the corresponding definitions of $\succ_D$,
one easily checks that one can have, for instance,
$(v,v') \succ_D (w, w)$ 
but, in contrast, $(w+\epsilon, w) \succ_D (v,v')$ for arbitrarily small $\epsilon>0$;
to this effect it suffices to assume $1>v>v'>0$ and to choose $w$ as follows depending on {\small$D$}: $w>v$ for {\small$D\,=\,\text{win}$}; $w<v'$ for {\small$D\,=\,\text{los}$}; $w=0$ for {\small$D\,=\,\text{ratio}$}.

\renewcommand\uplapar{\vskip-7mm\null}
\uplapar
\section{Examples}

\newcommand\vote[2]{\hss\ensep #2}
\newcommand\candidate[2]{\ensep\hbox to.7em{\hss\skl{#1}\hss}: #2}

\paragraph{3.1} As an example of a vote which involved truncated rankings, we look at an election which took place the 16th of February of 1652 in the Spanish royal household. This election is quoted in~\cite{ri}, but we use the more reliable and slightly different data which are given in~\cite[vol.\,2, p.\,263--264]{cv}. The office under election was that of ``aposentador mayor de palacio'', and the king was assessed
by six noblemen, who expressed respectively the following preferences about six candidates \skl{a}--\skl{f}:
\vote{Marqu\'es de Ari\c{c}a}{\skl{b}\better \skl{e}\better \skl{d}\better \skl{a}};
\vote{Conde de Barajas}{\skl{b}\better \skl{a}\better \skl{f}};
\vote{Conde de Montalb\'an}{\skl{a}\better \skl{f}\better \skl{b}\better \skl{d}};
\vote{Marqu\'es de Povar}{\skl{e}\better \skl{b}\better \skl{f}\better \skl{c}};
\vote{Conde de Pu\~nonrostro}{\skl{e}\better \skl{a}\better \skl{b}\better \skl{f}};\vote{Conde de Ysingui\'en}{\skl{b}\better \skl{d}\better \skl{a}\better \skl{f}}.

\medskip
The CLC computations are shown below.
Instead of relative scores, margins and turnouts, we will show their absolute counterparts, \ie without dividing them by the total number of votes.
This has the virtue of staying with small integer numbers.
Pairwise information will be given by means of square matrices;
as usual, the cell in row $x$ and column $y$ corresponds to the pair $xy$.
Since we consider only pairs $xy$ with $x\neq y$, we use the diagonal cells for specifying the simultaneous labelling of rows and columns by the members of $\ist$.
We must specify the labelling since at a certain stage we will rearrange the options in accordance with an admissible order.

We begin by forming the Llull matrix $(V_{xy})$ (step~0)
and deriving the corresponding indirect scores $(V^\ast_{xy})$ (step~1):
\begin{equation}
(V_{xy}) \,=\, 
\begin{small}
\begin{tabular}{|c|c|c|c|c|c|}
\hlinestrut
\labelcell{a}&2&5&3&3&5\\
\hlinestrut
4&\labelcell{b}&6&6&4&5\\
\hlinestrut
1&0&\labelcell{c}&1&0&0\\
\hlinestrut
2&0&3&\labelcell{d}&2&2\\
\hlinestrut
3&2&3&3&\labelcell{e}&3\\
\hlinestrut
1&1&5&4&3&\labelcell{f}\\
\hline
\end{tabular}
\end{small}
\,;
\qquad
(V^\ast_{xy}) \,=\, 
\begin{small}
\begin{tabular}{|c|c|c|c|c|c|}
\hlinestrut
\labelcell{a}&2&\bf5&\bf4&\bfgr3&\bf5\\
\hlinestrut
\bf4&\labelcell{b}&\bf6&\bf6&\bf4&\bf5\\
\hlinestrut
1&1&\labelcell{c}&1&1&1\\
\hlinestrut
2&2&\bf3&\labelcell{d}&2&2\\
\hlinestrut
\bfgr3&2&\bf3&\bf3&\labelcell{e}&\bfgr3\\
\hlinestrut
3&2&\bf5&\bf4&\bfgr3&\labelcell{f}\\
\hline
\end{tabular}
\end{small}
\,.
\end{equation}
In the matrix of the indirect scores we have visualized the indirect comparison relation
$\icr$ by marking in black bold face those pairs $xy$ that belong to it, \ie that satisfy $V^\ast_{xy} > V^\ast_{yx}$. Those that satisfy $V^\ast_{xy} = V^\ast_{yx}$ are shown in grey bold face. In accordance with (\ref{eq:copeland}), the Copeland rank of each option is then easily worked out as the number of black bold faces in the corresponding column \,plus half the number of grey bold faces in it \,plus one. The resulting values are shown next:
\begin{equation}
(\rank{x}) \,=\, 
\begin{small}
\begin{tabular}{|c|c|c|c|c|c|}
\hlinestrut
\labelcell{a}&\labelcell{b}&\labelcell{c}&\labelcell{d}&\labelcell{e}&\labelcell{f}\\
\hlinestrut
2\pq12&1&6&5&3&3\pq12\\
\hline
\end{tabular}
\end{small}
\,.
\end{equation}
Arranging the options by non-decreasing values of $\rank{x}$ gives an admissible order.
In this case there is only one possibility, namely \skl{b}\better\skl{a}\better\skl{e}\better\skl{f}\better\skl{d}\better\skl{c}. This completes step~2. From now on, the options will be arranged in this order. We now take the indirect margins $M^\mast_{xy}=V^\ast_{xy} - V^\ast_{yx}$ to work out the superdiagonal intermediate projected margins $M^\sigma_{xx'}$ (step~3). Having arranged the options in the adopted admissible order, $M^\sigma_{xx'}$ is simply the minimum $M^\mast$ value in the rectangle that lies to the right of $x$ and above $x'$:
\begin{equation}
\hskip-1pt 
(M^\mast_{xy}) \,=\, 
\begin{small}
\begin{tabular}{|c|c|c|c|c|c|}
\hlinestrut
\labelcell{b}&2&2&3&4&5\\
\hlinestrut
\st&\labelcell{a}&0&2&2&4\\
\hlinestrut
\st&\st&\labelcell{e}&0&1&2\\
\hlinestrut
\st&\st&\st&\labelcell{f}&2&4\\
\hlinestrut
\st&\st&\st&\st&\labelcell{d}&2\\
\hlinestrut
\st&\st&\st&\st&\st&\labelcell{c}\\
\hline
\end{tabular}
\end{small}
\,;
\hskip-1.5pt 
\qquad
(M^\sigma_{xx'}) \,=\, 
\begin{small}
\begin{tabular}{|c|c|c|c|c|c|}
\hlinestrut
\labelcell{b}&2&\st&\st&\st&\st\\
\hlinestrut
\st&\labelcell{a}&0&\st&\st&\st\\
\hlinestrut
\st&\st&\labelcell{e}&0&\st&\st\\
\hlinestrut
\st&\st&\st&\labelcell{f}&1&\st\\
\hlinestrut
\st&\st&\st&\st&\labelcell{d}&2\\
\hlinestrut
\st&\st&\st&\st&\st&\labelcell{c}\\
\hline
\end{tabular}
\end{small}
\,.
\end{equation}
We now take the original turnouts $T_{xy}$ to determine the intermediate projected ones $T^\sigma_{xy}$ (step~4). This involves also the superdiagonal intermediate projected margins $M^\sigma_{xx'}$ and requires solving (the absolute counterpart of) the problem of minimizing (\ref{eq:phi}) under the constraints (\ref{eq:symmetry}--\ref{eq:increment}). The results are shown below:
\begin{equation}
(T_{xy}) \,=\, 
\begin{small}
\begin{tabular}{|c|c|c|c|c|c|}
\hlinestrut
\labelcell{b}&6&6&6&6&6\\
\hlinestrut
\st&\labelcell{a}&6&6&5&6\\
\hlinestrut
\st&\st&\labelcell{e}&6&5&3\\
\hlinestrut
\st&\st&\st&\labelcell{f}&6&5\\
\hlinestrut
\st&\st&\st&\st&\labelcell{d}&4\\
\hlinestrut
\st&\st&\st&\st&\st&\labelcell{c}\\
\hline
\end{tabular}
\end{small}
\,;
\qquad
(T^\sigma_{xy}) \,=\, 
\begin{small}
\begin{tabular}{|c|c|c|c|c|c|}
\hlinestrut
\labelcell{b}&6&6&6&6&6\\
\hlinestrut
\st&\labelcell{a}&6&6&5\pq13&4\pq23\\
\hlinestrut
\st&\st&\labelcell{e}&6&5\pq13&4\pq23\\
\hlinestrut
\st&\st&\st&\labelcell{f}&5\pq13&4\pq23\\
\hlinestrut
\st&\st&\st&\st&\labelcell{d}&4\\
\hlinestrut
\st&\st&\st&\st&\st&\labelcell{c}\\
\hline
\end{tabular}
\end{small}
\,.
\end{equation}
Although the last step involved all of the intermediate projected turnouts,
the next one uses only the superdiagonal ones $T^\sigma_{xx'}$. These numbers together with the superdiagonal intermediate margins $M^\sigma_{xx'}$ determine the intervals
$\varGamma_{xx'}=((T^\sigma_{x'x}-M^\sigma_{xx'})/2,(T^\sigma_{x'x}+M^\sigma_{xx'})/2)=(V^\pi_{x'x},V^\pi_{xx'})$, whose unions as in~(\ref{eq:gammaxy}) give all the other intervals $\varGamma_{xy}=(V^\pi_{yx},V^\pi_{xy})$ (step~5):
\begin{equation}
(V^\pi_{xy}) \,=\, 
\begin{small}
\begin{tabular}{|c|c|c|c|c|c|}
\hlinestrut
\labelcell{b}&4&4&4&4&4\\
\hlinestrut
2&\labelcell{a}&3&3&3\pq16&3\pq16\\
\hlinestrut
2&3&\labelcell{e}&3&3\pq16&3\pq16\\
\hlinestrut
2&3&3&\labelcell{f}&3\pq16&3\pq16\\
\hlinestrut
2&2\pq16&2\pq16&2\pq16&\labelcell{d}&3\\
\hlinestrut
1&1&1&1&1&\labelcell{c}\\
\hline
\end{tabular}
\end{small}
\,.
\end{equation}
Finally, the rank-like rates are obtained from the projected scores by means of formula (\ref{eq:rrates}) with $\psc_{xy} = V^\pi_{xy}/V$ (step~6):
\begin{equation}
(\rlr{x}) \,=\, 
\begin{small}
\begin{tabular}{|c|c|c|c|c|c|}
\hlinestrut
\labelcell{b}&\labelcell{a}&\labelcell{e}&\labelcell{f}&\labelcell{d}&\labelcell{c}\\
\hlinestrut
2.6667&
3.6111&
3.6111&
3.6111&
4.0833&
5.1667\\
\hline
\end{tabular}
\end{small}
\,.
\end{equation}

\medskip
According to these results, the office should have been given to candidate~\skl{b}, who is also the winner by most other methods.
In the CLC method, this candidate is followed by three runners-up tied to each other, namely candidates \skl{a},\skl{e} and \skl{f}. In spite of the clear advantage of candidate~\skl{b}, the king appointed candidate \skl{f}, which was the celebrated painter Diego Vel\'azquez.

\paragraph{3.2} As a second example of an election involving truncated rankings we take the Debian Project leader election, which is using the qualitative method of Markus Schulze since 2003.
So far, the winners of these elections have been clear enough.
However, a quantitative measure of this clearness was lacking.
In the~following we consider the 2006 election, which had a participation of $V=421$ actual voters out of a total population of $972$ members. The individual votes  were taken from \texttt{\small http://www.debian.org/vote/2006/vote\underl002}.

\medskip
The next tables show the Llull matrix of that election and the resulting rank-like rates:

\begin{equation}
(V_{xy}) \,=\, 
\begin{small}
\begin{tabular}{|c|c|c|c|c|c|c|c|}
\hlinestrut
\labelcell{1}&\bf321&144&159\pq12&\bfgr{193\pq12}&\bf347\pq12&\bf246&\bf320\\
\hlinestrut
51&\labelcell{2}&42&53&50&\bf262&65&163\\
\hlinestrut
\bf251&\bf340&\labelcell{3}&198\pq12&\bf253&\bf362&\bf300&\bf345\\
\hlinestrut
\bf245\pq12&\bf341&\bf204\pq12&\labelcell{4}&\bf256&\bf371\pq12&\bf291\pq12&\bf339\pq12\\
\hlinestrut
\bfgr{193\pq12}&\bf325&144&149&\labelcell{5}&\bf357&\bf254&\bf321\pq12\\
\hlinestrut
26\pq12&77&24&22\pq12&21&\labelcell{6}&30&74\pq12\\
\hlinestrut
137&\bf292&90&109\pq12&131&\bf330&\labelcell{7}&\bf296\\
\hlinestrut
76&\bf207&54&71\pq12&75\pq12&\bf302\pq12&89&\labelcell{8}\\
\hline
\end{tabular}
\end{small}
\,,
\end{equation}

\begin{equation}
(\rlr{x}) \,=\, 
\begin{small}
\begin{tabular}{|c|c|c|c|c|c|c|c|}
\hlinestrut
\labelcell{1}&\labelcell{2}&\labelcell{3}&\labelcell{4}&\labelcell{5}&\labelcell{6}&\labelcell{7}&\labelcell{8}\\
\hlinestrut
\hbox to2.25em{\hss4.1105\hss}&
\hbox to2.25em{\hss5.9145\hss}&
\hbox to2.25em{\hss3.6926\hss}&
\hbox to2.25em{\hss3.6784\hss}&
\hbox to2.25em{\hss4.1105\hss}&
\hbox to2.25em{\hss6.7197\hss}&
\hbox to2.25em{\hss4.5720\hss}&
\hbox to2.25em{\hss5.8100\hss}\\
\hline
\end{tabular}
\end{small}
\,.
\end{equation}

%
%

\bigskip\noindent
In this case, the CLC results are in full agreement with the Copeland ranks of the original Llull matrix. In particular, both of them give an exact tie between candidates~\skl{1} and \skl{5}. Even so,
the CLC rates yield a quantitative information which is not present in the Copeland ranks.


For the computation of the rates we have taken $V=421$ (the actual number of votes) instead of $V=972$ (the number of people with the right to vote).
This is especially justified in Debian elections since they systematically include ``none~of~the~above'' as one of the alternatives, so it is reasonable to~interpret that abstention does not have a critical character. In~this case, ``none~of~the~above'' was alternative~\skl{8}, which 
obtained a~better result than two of the real candidates.

\paragraph{3.3} Finally, we look at an example of approval voting. Specifically, we consider the 2006 Public Choice Society election~\cite{pcs06}. Besides an approval vote, here the voters were also asked for a preferential vote ``in the spirit of research on public choice''. However,  here we will limit ourselves to the approval vote, which was the official one.
The vote had a participation of $V = 37$ voters, most of which approved more than one candidate.

The actual ballots are listed in the following table,\footnote{We are grateful to Prof.\ Steven J.\ Brams, who was the president of the Public Choice Society when that election took place, for his kind permission to reproduce these data.} where we give not only the approval voting data but also the the expressed preferences. The approved candidates are the ones which lie at the left of the slash.

\newcommand\fiav{\,/\,}

\bigskip
\vbox{%
\leavevmode
\hbox to125mm{\hrulefill}\par

\bgroup\small
\centerline{%
\hfil
\vtop{\noindent
\halign{\strut#\hfill\cr
\skl{A}\better \skl{B}\fiav\cr 
\skl{A}\better \skl{C}\better \skl{B}\fiav\cr 
\skl{D}\fiav \skl{A}\better \skl{B}\better \skl{E}\better \skl{C}\cr 
\skl{B}\better \skl{A}\fiav \skl{D}\better \skl{C}\better \skl{E}\cr 
\skl{D}\better \skl{A}\better \skl{B}\better \skl{C}\fiav \skl{E}\cr 
\skl{C}\better \skl{B}\better \skl{A}\fiav\cr 
\skl{E}\fiav \skl{D}\cr 
\skl{C}\better \skl{A}\better \skl{B}\better \skl{E}\fiav\cr 
\skl{D}\better \skl{E}\fiav \skl{C}\better \skl{A}\better \skl{B}\cr 
\skl{E}\fiav\cr 
\skl{B}\better \skl{C}\fiav\cr 
\skl{D}\better \skl{C}\fiav \skl{B}\better \skl{E}\better \skl{A}\cr 
\skl{B}\fiav\cr 
}}
\hfil
\vtop{\noindent
\halign{\strut#\hfill\cr
\skl{A}\fiav\cr 
\skl{A}\fiav\cr 
\skl{D}\fiav \skl{A}\tied \skl{B}\tied \skl{C}\tied \skl{E}\cr 
\skl{A}\tied \skl{C}\fiav\cr 
\fiav\skl{B}\better \skl{E}\better \skl{A}\better \skl{D}\better \skl{C}\cr 
\skl{A}\tied \skl{B}\tied \skl{E}\fiav\cr 
\skl{A}\tied \skl{B}\tied \skl{C}\tied \skl{D}\tied \skl{E}\fiav\cr 
\skl{D}\better \skl{A}\better \skl{B}\fiav\cr 
\skl{B}\better \skl{D}\better \skl{A}\fiav \skl{C}\better \skl{E}\cr 
\skl{A}\fiav \skl{B}\better \skl{E}\better \skl{C}\better \skl{D}\cr 
\skl{D}\fiav\cr 
\skl{A}\tied \skl{C}\better \skl{B}\fiav \skl{D}\better \skl{E}\cr 
\skl{A}\fiav \skl{D}\better \skl{B}\better \skl{C}\better \skl{E}\cr 
}}
\hfil
\vtop{\noindent
\halign{\strut#\hfill\cr
\skl{C}\fiav \skl{B}\better \skl{D}\better \skl{A}\better \skl{E}\cr 
\skl{C}\fiav\cr 
\skl{D}\tied \skl{E}\fiav \skl{A}\better \skl{B}\tied \skl{C}\cr 
\skl{B}\fiav \skl{C}\better \skl{A}\better \skl{D}\better \skl{E}\cr 
\skl{D}\better \skl{C}\better \skl{E}\fiav\cr 
\skl{C}\fiav \skl{A}\better \skl{B}\tied \skl{D}\tied \skl{E}\cr 
\skl{C}\fiav\cr 
\skl{B}\better \skl{D}\fiav \skl{E}\better \skl{C}\better \skl{A}\cr 
\skl{B}\better \skl{C}\fiav \skl{A}\better \skl{E}\better \skl{D}\cr 
\skl{D}\better \skl{A}\better \skl{C}\better \skl{B}\fiav\cr 
\skl{D}\better \skl{E}\fiav \skl{A}\better \skl{B}\cr 
}}
\hfil}
\egroup

\leavevmode
\hbox to125mm{\hrulefill}
}

\smallskip
\medskip
The next table gives the number of received approvals~$A_x$ together with the rank-like rates resulting from the four variants of the CLC method (the superindices $1,2,3,4$ indicate respectively the main, codual, balanced and margin-based variants).

\begin{equation*}
\begin{small}
\raise2ex\hbox{%
\begin{tabular}{|c|c|c|c|c|c|}
\hlinestrut\rule[-1.1ex]{0pt}{1ex}
$x$&\labelcell{A}&\labelcell{B}&\labelcell{C}&\labelcell{D}&\labelcell{E}\\
\hlinestrut\rule[-1.1ex]{0pt}{1ex}
$A_x$&17&16&17&14&9\\
\hlinestrut\rule[-1.1ex]{0pt}{1ex}
$\rlr{x}^1$&3.6014&3.6486&3.6149&3.7720&4.1689\\
\hlinestrut\rule[-1.1ex]{0pt}{1ex}
$\rlr{x}^2$&3.6081&3.6486&3.6081&3.7568&4.2162\\
\hlinestrut\rule[-1.1ex]{0pt}{1ex}
$\rlr{x}^3$&3.6081&3.6486&3.6081&3.7703&4.1622\\
\hlinestrut\rule[-1.1ex]{0pt}{1ex}
$\rlr{x}^4$&2.8919&2.9324&2.8919&3.0135&3.2703\\
\hline
\end{tabular}%
}
\end{small}
\end{equation*}

%



\medskip
\noindent
As one can see, the approval scores result in a tie for the first place between candidates~\skl{A} and \skl{C}, which are followed at a minimum distance by candidate~\skl{B} and then by candidates~\skl{D} and \skl{E}.
\ensep
Exactly the same ranking is found in the results of the codual, balanced and margin-based variants of the CLC~method, but not in those of the main variant, which discriminates between candidates~\skl{A} and \skl{C}, giving the victory to~\skl{A}.
\ensep
In~\secpar{\ref{sec-approval}} we will see that the ranking given by the margin-based variant is always in full agreement with that given by the approval score.


\section{The projection. Well-definedness and structural properties}
\label{sec-projection}

As it was seen in \secpar{2}, the CLC method involves in a crucial way a projection $(v_{xy})\mapsto(\psc_{xy})$ of the original Llull matrix onto a special set of such matrices.
\ensep
In this section we will see that this projection is well defined and we will look at certain structural properties of the resulting matrix.

\medskip
The projection comprises steps~1--5 of the CLC procedure.
In order to ensure that it is well defined, one must check the following points:\linebreak 
\ensep
(a)~The existence (and effective production) of an admissible order~$\xi$. Since it involves only steps~1--2, before any divergence of the procedure from the complete case, this point requires no other considerations than those made in~\cite{crc}.\linebreak 
\ensep
(b)~The existence and uniqueness of a minimizer of (\ref{eq:phi})
under the conditions~(\ref{eq:symmetry}--\ref{eq:increment}). This is  a consequence of the fact that the set $\pptoset$ defined by these conditions is a closed convex set \cite[ch.\,I, \secpar{2}]{ks}. In this connection, one can say that $(\ppto_{xy})$ is the orthogonal projection of $(t_{xy})$ onto the convex set~$\pptoset$.
\ensep
And (c)~that the final results are independent of the admissible order~$\xi$ when there are several possibilities for it. This will be dealt with by Theorem~\ref{st:independenceOfXi} below, whose proof is rather long, but not difficult.

\medskip
Before embarking on that theorem, however, we will look at certain properties of the intervals~$\gamma_{xy}$.
Besides being used in that theorem, these properties will be seen later on to be at the core of the structure of the projected Llull matrix. We will use the following notation:
$|\gamma|$~means the length of an interval,
and $\centre\gamma$~means its barycentre, or centroid,
\ie the~number $(a+b)/2$ \,if $\gamma=[a,b]$.

\smallskip
\begin{lemma}\hskip.5em
\label{st:intervals}
The sets $\gamma_{xy}$ have the following properties for~$x\rxi y\rxi z$:

\iim{a} $\gamma_{xy}$ is a closed interval.

\iim{b} $\gamma_{xy} \,\,\sbseteq\,\, [0,1]$.

\iim{c} $\gamma_{xz} \,\,=\,\, \gamma_{xy} \,\cup\, \gamma_{yz}$.

\iim{d} $\gamma_{xy} \,\cap\, \gamma_{yz} \,\,\neq\,\, \emptyset$.

\iim{e} $|\gamma_{xz}| \,\,\ge\,\, \max\,(\,|\gamma_{xy}|\,,\, |\gamma_{yz}|\,)$.

\iim{f} $\centre\gamma_{xy} \,\,\ge\,\, \centre\gamma_{xz} \,\,\ge\,\,\centre\gamma_{yz}$.

\iim{g} $\centre\gamma_{xy} \,-\, |\gamma_{yz}|/2 \,\,\le\,\, \centre\gamma_{xz} \,\,\le\,\, \centre\gamma_{yz} \,+\, |\gamma_{xy}|/2$.
\end{lemma}
\begin{proof}\hskip.5em
Let us begin by noticing that the superdiagonal intermediate turnouts and margins are ensured to satisfy the following inequalities:
\begin{gather}
\label{eq:bounded0}
0 \,\le\, \ppmg_{xx'} \,\le\, \ppto_{xx'} \,\le\, 1
\\[2.5pt]
\label{eq:overlap}
0 \,\le\, \ppto_{xx'} - \ppto_{x'x''} \,\le\, \ppmg_{xx'} + \ppmg_{x'x''}.
\end{gather}
The inequalities of~(\ref{eq:bounded0}) are those of (\ref{eq:bounded}), with $\ppto_{xx'}$ substituted for $\gto_{xx'}$, plus the fact that $\ppmg_{xx'}\ge0$.
Those of~(\ref{eq:overlap}) are the result of adding up (\ref{eq:increment}) with $x$ and $z$ replaced respectively by $x'$ and $x$ \,plus\, the same inequality with $z$ replaced by $x''$, and using the symmetric character of the turnouts.

\halfsmallskip
From (\ref{eq:bounded0}) it follows that
$0\le(\ppto_{xx'}-\ppmg_{xx'})/2\le(\ppto_{xx'}+\ppmg_{xx'})/2\le1$.
So, every $\gamma_{xx'}$ is an interval (possibly reduced to one point) and this interval is contained in $[0,1]$.
Also, the inequalities of~(\ref{eq:overlap}) ensure on the one hand that $\centre\gamma_{xx'}\ge\centre\gamma_{x'x''}$, and on the other hand that the intervals $\gamma_{xx'}$ and $\gamma_{x'x''}$ overlap each other.
\ensep
In the following we will see that these facts about the elementary intervals $\gamma_{xx'}$ entail the stated properties of the sets $\gamma_{xy}$ defined by~(\ref{eq:gammaxy}).

\newcommand\ppart[1]{\halfsmallskip\leavevmode
\hbox to4.25em{Part~(#1).\hss}\ignorespaces}

\ppart{a} This is an obvious consequence of the fact that $\gamma_{pp'}$ and $\gamma_{p'p''}$ overlap each other.

\ppart{b} This follows from the fact that $\gamma_{pp'}\sbseteq[0,1]$.

\ppart{c} This is a consequence of the associative property enjoyed by the set-union operation.

\ppart{d} This is again an obvious consequence of the fact that $\gamma_{pp'}$ and $\gamma_{p'p''}$ overlap each other (take $p'=y$).

\ppart{e} This follows from~(c) because $\intg\sbseteq\inth$ implies $|\intg|\le|\inth|$.

\ppart{f} This follows from the inequality $\centre\gamma_{pp'}\ge\centre\gamma_{p'p''}$ because of the following general fact:
\ensep
If~$\intg$ and~$\inth$ are two overlapping intervals with $\centre\intg\ge\centre\inth$ then 
$\centre\intg \ge \pcentre{\intg\cup\inth} \ge \centre\inth$. This follows immediately from the definition of the barycentre.

\ppart{g} This follows from (c) and (d) because of the following general fact:\ensep
If~$\intg$ and~$\inth$ are two overlapping intervals, then 
$\centre\intg - |\inth|/2 \le \pcentre{\intg\cup\inth} \le \centre\intg + |\inth|/2$. Again, this follows easily from the definitions.
%
\end{proof}

\smallskip
\begin{theorem}\hskip.5em
\label{st:independenceOfXi}
The projected scores do not depend on the admissible order~$\xi$ used for their calculation,
\ie the value of\, $\psc_{xy}$ is independent of\, $\xi$ \,for every $xy\in\tie$.
On the other hand, the matrix of the projected scores
in an admissible order~$\xi$
is also independent of\, $\xi$;
\ie if $x_i$ denotes the element of rank~$i$ in~$\xi$,
the value of $\psc_{x_ix_j}$ is independent of\, $\xi$
\,for every pair of indices~$i,j$.
\end{theorem}

\noindent
\textit{Remark}.\hskip.5em
The two statements say different things since the identity of $x_i$ and $x_j$ may depend on the admissible order~$\xi$.

\begin{proof}\hskip.5em
For the purposes of this proof it becomes necessary to change our set-up in a certain way. In fact, until now the intermediate objects $\ppmg_{xy}$, $\ppto_{xy}$ and $\gamma_{xy}$ were considered only for~$x\rxi y$, \ie $xy\in\xi$. However, since we have to deal with changing the admissible order $\xi$, here we will allow their argument $xy$ to be any pair (of different elements), no matter whether it belongs to $\xi$ or not. In~this connection, we will certainly put $\ppmg_{yx}=-\ppmg_{xy}$ and $\ppto_{yx}=\ppto_{xy}$. On~the other hand, concerning $\gamma_{xy}$ and $\gamma_{yx}$, we will proceed in the following way: if $\gamma_{xy}=[a,b]$ then $\gamma_{yx}=[b,a]$. So,~generally speaking the $\gamma_{xy}$ are here ``oriented intervals'', \ie ordered pairs of real numbers. However, $\gamma_{xy}$ will always be ``positively oriented'' when $xy$ belongs to an admissible order (but it will be reduced to a point whenever there is another admissible order which includes~$yx$). In~particular, the $\gamma_{pp'}$ which are combined in (\ref{eq:gammaxy}) are always positively oriented intervals; so, the union operation performed in that equation can always be understood in the usual sense.
\ensep
In the following, $\gamma\rev$ denotes the oriented interval ``reverse'' to $\gamma$, \ie $\gamma\rev=[b,a]$ if $\gamma=[a,b]$.

So, let us consider the effect of replacing $\xi$ by another admissible order~$\xibis$. In the following, the tilde is systematically used to distinguish between hom\-olo\-gous objects which are associated respectively with $\xi$~and~$\xibis$; in particular, such a notation will be used in connection with the labels of the equations which are formulated in terms of
the assumed admissible order.

With this terminology, we will prove the two following equalities. First,
\begin{equation}
\hskip.75em\gamma_{xy} \,=\, \gammabis_{xy},\hskip.75em\qquad
\hbox to60mm{for any pair\, $xy\,\ (x\neq y)$,\hfil}
\label{eq:ginvariance}
\end{equation}
where $\gamma_{xy}$ are the intervals produced by (\ref{eq:psi}--\ref{eq:gammaxy}) together with the operation $\gamma_{yx}=\gpxy$,
and $\gammabis_{xy}$ are those produced by (\tref{eq:psi}--\tref{eq:gammaxy}) together with the operation $\gammabis_{yx}=\gbispxy$.
\ensep
Secondly, we will see also that
\begin{equation}
\gamma_{x_ix_j} \,=\, \gammabis_{\tilde x_i\tilde x_j},\qquad
\hbox to60mm{for any pair of indices\, $ij\,\ (i\neq j),$\hfil}
\label{eq:gijinvariance}
\end{equation}
where $x_i$ denotes the element of rank $i$ in $\xi$, and analogously for~$\tilde x_i$ in $\xibis$.
\ensep
These equalities contain the statements of the theorem since the projected scores are nothing else than the end points of the intervals~$\gamma_{xy}$.

\medskip
Now, by a well-known result, proved for instance in \cite{fh}, 
it suffices to deal
with the case of two admissible orders $\xi$~and~$\xibis$ which
differ from each other by one inversion only.
So, we will assume that there are two elements $a$ and $b$ such that the only difference between $\xi$~and $\xibis$ is that $\xi$ contains $ab$ whereas $\xibis$~contains~$ba$.
According to the definition of an admissible order, this implies that
$\img_{ab} = \img_{ba} = 0$.

\newcommand\Ant{P}
\newcommand\ant{p}
\newcommand\Pos{Q}
\newcommand\pos{q}

In order to control the effect of the differences between $\xi$ and $\xibis$, we will make use of the following notation:
$\ant$ will denote 
the immediate predecessor of~$a$ in $\xi$; in~this connection, any statement about $\ant$ will be understood to imply the assumption that 
the set of predecessors of $a$ in $\xi$ is not empty. Similarly, 
$\pos$ will denote 
the immediate successor of~$b$ in $\xi$; here too, any statement about $\pos$ will be understood to imply the assumption that 
the set of successors of $b$ in $\xi$ is not empty. So, $\xi$ and $\xibis$ contain respectively the paths $\ant ab\pos$ and $\ant ba\pos$.

\medskip
Let us look first at the superdiagonal intermediate projected margins~$\ppmg_{hh'}$. Since their definition is the same as in the complete case, one can invoke the same  arguments as in the proof of Theorem~6.2 of \cite{crc} to obtain the equality
\begin{equation}
\ppmg_{x_ix_{i+1}} \,=\, \ppmgbis_{\tilde x_i\tilde x_{i+1}},\qquad \hbox{for any $i=1,2,\dots N\!-\!1$.}
\label{eq:mequivariance}
\end{equation}
In more specific terms,
\begin{align}
\ppmg_{xx'}\! &\,=\, \ppmgbis_{xx'},\qquad \hbox{whenever \,$x\neq \ant,a,b$},
\label{eq:mxxp}
\\[2.5pt]
\ppmg_{\ant a} &\,=\, \ppmgbis_{\ant b},
\label{eq:mcacb}
\\[2.5pt]
\ppmg_{ab} &\,=\, \ppmgbis_{ba} \,=\, 0,
\label{eq:mabba}
\\[2.5pt]
\ppmg_{b\pos} &\,=\, \ppmgbis_{a\pos}.
\label{eq:mbdad}
\end{align}
In connection with equation~(\ref{eq:mxxp}) it should be clear that \emph{for~$x\neq \ant,a,b$ the immediate successor $x'$ is the same in both orders $\xi$ and $\xibis$}.

\medskip
Next we will see that the intermediate projected turnouts $\ppto_{xy}$ are invariant with respect to $\xi$:
\begin{equation}
\ppto_{xy} \,=\, \pptobis_{xy},\qquad\hbox{for any pair\, } xy\,\ (x\neq y),
\label{eq:tinvariance}
\end{equation}
where $\ppto_{xy}$ are the numbers produced by (\ref{eq:psi}) together with the symmetry $\ppto_{yx}=\ppto_{xy}$,
and $\pptobis_{xy}$ are those produced by (\tref{eq:psi}) together with the symmetry $\pptobis_{yx}=\pptobis_{xy}$.

We will prove (\ref{eq:tinvariance}) by seeing that the set $\pptoset$ determined by conditions (\ref{eq:symmetry},\ref{eq:bounded},\ref{eq:increment}) coincides exactly with the set $\pptosetbis$ determined by (\ref{eq:symmetry},\tref{eq:bounded},\tref{eq:increment}). In~other words, conditions (\ref{eq:bounded}--\ref{eq:increment}) are exactly equivalent to (\tref{eq:bounded}--\tref{eq:increment}) under condition (\ref{eq:symmetry}), which does not depend on~$\xi$.

In order to prove this equivalence we begin by noticing that condition~(\ref{eq:bounded}) coincides exactly with (\tref{eq:bounded}) when $x\neq \ant,a,b$. This is true because, on the one hand, $x'$ is then the same in both orders $\xi$ and $\xibis$, and, on the other hand, (\ref{eq:mxxp}) ensures that the right-hand sides have the same value. Similarly happens with conditions~(\ref{eq:increment}) and (\tref{eq:increment}) when $z\neq \ant,a,b$.
So, it remains to deal with conditions (\ref{eq:bounded}) and (\tref{eq:bounded}) for~$x=\ant,a,b$, and with conditions (\ref{eq:increment}) and (\tref{eq:increment}) for~$z=\ant,a,b$.\ensep
Now, on account of the symmetry~(\ref{eq:symmetry}), one easily sees that  condition (\ref{eq:bounded}) with $x=a$ is equivalent to (\tref{eq:bounded}) with $x=b$. In~fact, both of them reduce to $0 \le \gto_{ab} \le 1$ since $\ppmg_{ab} = \ppmgbis_{ba} = 0$, as it was   obtained in (\ref{eq:mabba}). This last equality ensures also the equivalence between condition~(\ref{eq:increment}) with~$z=a$ and condition~(\tref{eq:increment}) with~$z=b$. In this case both of them reduce to
\begin{equation}
\gto_{xa} \,=\, \gto_{xb}.
\label{eq:xaxb}
\end{equation}
This common equality plays a central role in
the equivalence between the remaining conditions.
\ensep
Thus, its combination with (\ref{eq:mbdad}) ensures the equivalence between  (\ref{eq:bounded}) with $x=b$ and (\tref{eq:bounded}) with $x=a$, as well as the equivalence between (\ref{eq:increment}) with $z=b$ and (\tref{eq:increment}) with $z=a$ when $x\neq a,b$. 
\ensep
On the other hand, its combination with (\ref{eq:mcacb}) ensures the equivalence between (\ref{eq:bounded}) and (\tref{eq:bounded}) when $x=\ant$, as well as the equivalence between (\ref{eq:increment}) and (\tref{eq:increment}) when $z=\ant$ and $x\neq a,b$. 
\ensep
Finally, we have the two following equivalences: (\ref{eq:increment})~with $z=\ant$ and $x=b$ is equivalent to (\tref{eq:increment}) with $z=\ant$ and $x=a$ because of the same equality  (\ref{eq:xaxb}) together with~(\ref{eq:mcacb}) and the symmetry~(\ref{eq:symmetry});
and similarly, (\ref{eq:increment}) with $z=b$ and $x=a$ is equivalent to (\tref{eq:increment}) with $z=a$ and $x=b$ because of (\ref{eq:xaxb}) together with~(\ref{eq:mbdad}) and~(\ref{eq:symmetry}).
\ensep
This completes the proof of~(\ref{eq:tinvariance}).

\vskip-3pt 
\medskip
Having seen that condition~(\ref{eq:xaxb}) is included in both~(\ref{eq:increment}) and~(\tref{eq:increment}), it follows that the intermediate projected turnouts satisfy
\begin{equation}
\ppto_{xa} \,=\, \ppto_{xb},\qquad \pptobis_{xa} \,=\, \pptobis_{xb}.
\end{equation}
By taking $x=\ant,\pos$ and using also~(\ref{eq:tinvariance}), it follows that
\begin{align}
\ppto_{xx'}\! &\,=\, \pptobis_{xx'},\qquad \hbox{whenever \,$x\neq \ant,a,b$},
\label{eq:txxp}
\\[2.5pt]
\ppto_{\ant a} &\,=\, \pptobis_{\ant b},
\label{eq:tcacb}
\\[2.5pt]
\ppto_{ab} &\,=\, \pptobis_{ba},
\label{eq:tabba}
\\[2.5pt]
\ppto_{b\pos} &\,=\, \pptobis_{a\pos}.
\label{eq:tbdad}
\end{align}
In other words, the superdiagonal intermediate turnouts satisfy
\begin{equation}
\ppto_{x_ix_{i+1}} \,=\, \pptobis_{\tilde x_i\tilde x_{i+1}},\qquad \hbox{for any $i=1,2,\dots N\!-\!1$.}
\label{eq:tequivariance}
\end{equation}
On account of the definition of $\gamma_{x_ix_{i+1}}$ and $\gammabis_{\tilde x_i\tilde x_{i+1}}$,
the combination of (\ref{eq:mequivariance}) and (\ref{eq:tequivariance}) results in
\begin{equation}
\gamma_{x_ix_{i+1}} \,=\, \gammabis_{\tilde x_i\tilde x_{i+1}},\qquad \hbox{for any $i=1,2,\dots N\!-\!1$,}
\label{eq:gequivariance}
\end{equation}
from which the union operation~(\ref{eq:gammaxy}) produces~(\ref{eq:gijinvariance}).

\medskip
Finally, let us see that (\ref{eq:ginvariance}) holds too. To this effect, we begin by noticing that (\ref{eq:mabba}) together with (\ref{eq:tabba}) are saying not only that $\gamma_{ab} = \gammabis_{ba}$ but also that this interval reduces to a point. As a consequence, we have
\begin{equation}
\gamma_{ba} \,=\, \gamma_{ab} \,=\, \gammabis_{ba} \,=\, \gammabis_{ab}.
\label{eq:4gabba}
\end{equation}
Let us consider now the equation $\gamma_{\ant a} = \gammabis_{\ant b}$, which is contained in (\ref{eq:gequivariance}).
Since $\gamma_{ab}$ reduces to a point, parts~(c) and (d) of Lemma~\ref{st:intervals} give
$\gamma_{\ant b} = \gamma_{\ant a}\cup\gamma_{ab} = \gamma_{\ant a}$.
Analogously, $\gammabis_{pa} = \gammabis_{\ant b}\cup\gammabis_{ba} = \gammabis_{\ant b}$. Altogether, this gives
\begin{equation}
\gamma_{\ant b} \,=\, \gamma_{\ant a} \,=\, \gammabis_{\ant b} \,=\, \gammabis_{\ant a}.
\label{eq:4gcacb}
\end{equation}
By means of an analogous argument, one obtains also that
\begin{equation}
\gamma_{a\pos} \,=\, \gamma_{b\pos} \,=\, \gammabis_{a\pos} \,=\, \gammabis_{b\pos}.
\label{eq:4gbdad}
\end{equation}
On the other hand, (\ref{eq:gequivariance}) ensures that
\begin{equation}
\gamma_{xx'} \,=\, \gammabis_{xx'},\qquad \hbox{whenever \,$x\neq \ant,a,b$}.
\label{eq:4gxxp}
\end{equation}
Finally, part~(c) of Lemma~\ref{st:intervals} allows to go from (\ref{eq:4gabba}--\ref{eq:4gxxp}) to the desired general equality~(\ref{eq:ginvariance}).
\end{proof}

\smallskip
\begin{theorem}\hskip.5em
\label{st:propertiesOfProjection}
The projected scores and their asssociated margins and turn\-outs satisfy the following properties with respect to any admissible order~$\xi$:

\halfsmallskip\noindent
\textup{(a)}~The following inequalities hold for $x\rxi y\rxi z$:
\begin{gather}
\label{eq:pvxyinequality}
\psc_{xy} \,\ge\, \psc_{yx},\quad\text{\ie } \pmg_{xy} \,\ge\, 0,
\\[2.5pt]
\label{eq:pvequaltomax}
\psc_{xz} \,\,=\,\, \max\,(\psc_{xy},\psc_{yz}),
\\[2.5pt]
\label{eq:pvequaltomin}
\psc_{zx} \,\,=\,\, \min\,(\psc_{zy},\psc_{yx}),
\\[2.5pt]
\label{eq:pmlesthanmax}
\pmg_{xz} \,\,\le\,\, \pmg_{xy} \,+\, \pmg_{yz},
\\[2.5pt]
\label{eq:pttm}
\pto_{xz} \,-\, \pto_{yz} \,\le\, \pmg_{xy},\quad
\pto_{xy} \,-\, \pto_{xz} \,\le\, \pmg_{yz}.
\end{gather}

\halfsmallskip\noindent
\textup{(b)}~The following inequalities hold for $x\rxi y$ and $z\not\in\{x,y\}$:
\begin{alignat}{2}
\label{eq:pvinequalities}
\psc_{xz} \,&\ge\, \psc_{yz},
&\qquad
\psc_{zx} \,&\le\, \psc_{zy},
\\[2.5pt]
\label{eq:pminequalities}
\pmg_{xz} \,&\ge\, \pmg_{yz},
&\qquad
\pmg_{zx} \,&\le\, \pmg_{zy},
\\[2.5pt]
\label{eq:ptinequalities}
\pto_{xz} \,&\ge\, \pto_{yz},
&\qquad
\pto_{zx} \,&\ge\, \pto_{zy},
\end{alignat}

\halfsmallskip\noindent
\textup{(c)} If~$\psc_{xy}=\psc_{yx}$, or equivalently~$\pmg_{xy}=0$,  then \textup{(\ref{eq:pvinequalities}--\ref{eq:ptinequalities})} are satisfied all of them with an equality sign.

\halfsmallskip\noindent
\textup{(d)}~The absolute projected margins $\apm_{xy}=|\pmg_{xy}|$ satisfy the triangular 
inequality \,$\apm_{xz} \le \apm_{xy} + \apm_{yz}$\, for any $x,y,z$.
\end{theorem}

\begin{proof}\hskip.5em
We will see that these properties derive from those satisfied by the intervals $\gamma$, which are collected in Lemma~\ref{st:intervals}.
For the derivation one has to bear in mind that $\psc_{xy}$ and $\psc_{yx}$ are respectively the right and left end points of the interval $\gamma_{xy}$, and that $\pmg_{xy} = -\pmg_{yx}$ and $\pto_{xy} = \pto_{yx}$ are respectively the width and twice the barycentre of $\gamma_{xy}$.

\smallskip
Part~(a).\ensep
First of all, (\ref{eq:pvxyinequality}) holds as soon as $\gamma_{xy}$ is an interval, as it is ensured by part~(a) of Lemma~\ref{st:intervals}. 
\ensep
On the other hand, (\ref{eq:pvequaltomax}--\ref{eq:pvequaltomin}) are nothing else than a paraphrase of Lemma~\ref{st:intervals}.(c), in the same way as (\ref{eq:pvequaltomax}--\ref{eq:pmlesthanmax}) is a paraphrase of Lemma~\ref{st:intervals}.(e).
\ensep
Finally, the inequalities of (\ref{eq:pttm}) are those of Lemma~\ref{st:intervals}.(g).

\smallskip
Part~(b).\ensep
Let us begin by noticing that (\ref{eq:pminequalities}) will be an immediate consequence of (\ref{eq:pvinequalities}), since $\pmg_{xz} = \psc_{xz} - \psc_{zx}$ and $\pmg_{yz} = \psc_{yz} - \psc_{zy}$.
On the other hand, (\ref{eq:ptinequalities}.2) is equivalent to (\ref{eq:ptinequalities}.1). This equivalence  holds because the turnouts are symmetric. So, it remains to prove the inequalities (\ref{eq:pvinequalities}) and either (\ref{eq:ptinequalities}.1) or (\ref{eq:ptinequalities}.2). In order to prove them we will distinguish three cases, namely:\ensep
(i)~$x\rxi y\rxi z$;\ensep
(ii)~$z\rxi x\rxi y$;\ensep
(iii)~$x\rxi z\rxi y$.

Case~(i)\,: By part~(c) of Lemma~\ref{st:intervals}, in this case we
have $\gamma_{xz}\spseteq\gamma_{yz}$. This immediately implies
(\ref{eq:pvinequalities}) because $[a,b]\spseteq[c,d\,]$ is
equivalent to saying that $b\ge d$ and $a\le c$. On the other hand,
the inequality~(\ref{eq:ptinequalities}.1) is contained in part~(f) of
Lemma~\ref{st:intervals}. Case~(ii)\, is analogous to case~(i).

Case~(iii)\,: In this case, (\ref{eq:pvinequalities}) follows from part~(d) of Lemma~\ref{st:intervals} since $[a,b]\cap[c,d\hskip2pt]\neq \emptyset$ is equivalent to saying that $b\ge c$ and $a\le d$. On the other hand,  (\ref{eq:ptinequalities}.1) is still contained in part~(f) of Lemma~\ref{st:intervals} (because of the symmetric character of the turnouts).

\smallskip
Part~(c).\ensep
The hypothesis that $\psc_{xy}=\psc_{yx}$ is equivalent to saying that $\gamma_{xy}$ reduces to a point, \ie $\gamma_{xy}=[v,v]$ for some $v$.
We will distinguish the same three cases as in part~(a).

Case~(i)\,: On account of the overlapping property $\gamma_{xy}\cap\gamma_{yz}\neq \emptyset$ (part~(d) of Lemma~\ref{st:intervals}), the one-point interval $\gamma_{xy}=[v,v]$ must be contained in $\gamma_{yz}$. So,~$\gamma_{xz} = \gamma_{xy}\cup\gamma_{yz} = \gamma_{yz}$ (where we used part~(c) of Lemma~\ref{st:intervals}).
Case~(ii)\, is again analogous to case~(i).

Case~(iii)\,: By part~(c) of Lemma~\ref{st:intervals} (with $y$ and $z$ interchanged with each other), the fact that $\gamma_{xy}$ reduces to the one-point interval $[v,v]$ implies that  both $\gamma_{xz}$ and $\gamma_{zy}$ reduce also to this one-point interval.

\smallskip
Part~(d).\ensep
It suffices to consider the particular ordering \,$x\rxi y\rxi z$\, and check that each of the three numbers $\apm_{xy}=\apm_{yx}=\pmg_{xy}$, $\apm_{yz}=\apm_{zy}=\pmg_{yz}$, $\apm_{xz}=\apm_{zx}=\pmg_{xz}$ is less than the sum of the other two.
\ensep
This is so since we know that $\pmg_{xz} \le \pmg_{xy}+\pmg_{yz}$, by (\ref{eq:pmlesthanmax}), and also that $\pmg_{yz} \le \pmg_{xz}$ and $\pmg_{xy} \le \pmg_{xz}$, by (\ref{eq:pminequalities}), or more directly, by Lemma~\ref{st:intervals}.(e).
\end{proof}

\bigskip
In \secpar{2.2.2} we have seen that in the case of single-choice voting the projected scores coincide with the original ones.
In that connection, one has the following general result:

\smallskip
\begin{proposition}\hskip.5em
\label{st:nochange}
Assume that there exists a total order $\xi$ such that
the original scores and their associated margins and turnouts satisfy the following conditions:
\begin{alignat}{2}
\label{eq:vxyinequality}
&v_{xy} \ge v_{yx},\text{ i.e.\ } m_{xy} \ge 0,\qquad &&\hbox{whenever $x\rxi y$},
\\[2.5pt]
\label{eq:vequaltomax}
&v_{xz} \,\,=\,\, \max\,(v_{xy},v_{yz}),\qquad &&\hbox{whenever $x\rxi y\rxi z$},
\\[2.5pt]
\label{eq:vequaltomin}
&v_{zx} \,\,=\,\, \min\,(v_{zy},v_{yx}),\qquad &&\hbox{whenever $x\rxi y\rxi z$},
\\[2.5pt]
\label{eq:ttm}
&0 \,\le\, t_{xz} - t_{x'z} \,\le\, m_{xx'},\qquad &&\hbox{whenever $z\not\in\{x,x'\}$}.
\end{alignat}
In that case, the projected scores coincide with the original ones.
\end{proposition}


\begin{proof}\hskip.5em
Let us begin by noticing that condition~(\ref{eq:ttm}) implies, as in the proof of Lemma~\ref{st:intervals}, that the intervals $[v_{x'x},v_{xx'}]$ and $[v_{x''x'},v_{x'x''}]$ overlap each other. In other words, one has the following inequalities ``accross the diagonal'': $v_{x'x}\le v_{x'x''}$ and $v_{x''x'}\le v_{xx'}$.
\ensep
Now, (\ref{eq:vequaltomax}--\ref{eq:vequaltomin}) together with the preceding inequalities imply that
\begin{equation}
\label{eq:originalvinequalities}
v_{xz} \,\ge\, v_{yz},
\quad
v_{zx} \,\le\, v_{zy},
\qquad
\hbox{whenever $x\rxi y$ and $z\not\in\{x,y\}$.}
\end{equation}
From these facts, one can derive that 
\begin{equation}
\label{eq:minInequalityIndirect}
v_{xz} \,\ge\, \min\,(v_{xy}, v_{yz}),\qquad\hbox{for any $x,y,z$.}
\end{equation}
In fact, for $z\rxi y\rxi x$ this inequality is guaranteed by~(\ref{eq:vequaltomin}) (with $x$ and $z$ interchanged with each other), whereas for any other ordering of $x,y,z$ one easily arrives at (\ref{eq:minInequalityIndirect}) as a consequence of~(\ref{eq:originalvinequalities}).

Now, according to Lemma~4.2 of \cite{crc}, (\ref{eq:minInequalityIndirect}) implies that the indirect scores coincide with the original ones: $\isc_{xy}=v_{xy}$, which entails that $\img_{xy}=m_{xy}$.\ensep
In particular, $\xi$ is ensured to be an admissible order.
\ensep

Proceeding with the CLC algorithm, one easily checks that $\ppmg_{xy}=m_{xy}$ (because the pattern of growth~(\ref{eq:originalvinequalities}) gets transmitted from the scores to the margins) and $\ppto_{xy}=t_{xy}$ (since the turnouts are assumed to satisfy (\ref{eq:ttm}) and one certainly has $m_{xy}\le t_{xy}\le1$). As a consequence of these facts, (\ref{eq:pscxxp}) results in $\psc_{xx'}=v_{xx'}$ and $\psc_{x'x}=v_{x'x}$, from which (\ref{eq:pscxy}) and (\ref{eq:vequaltomax}--\ref{eq:vequaltomin}) lead to conclude that $\psc_{xy}=v_{xy}$ for any pair $xy$.
%
\end{proof}

\smallskip
Since conditions (\ref{eq:vxyinequality}--\ref{eq:ttm}) of Proposition~\ref{st:nochange} are included among the properties of the projected Llull matrix according to Theorem~\ref{st:propertiesOfProjection}, one can conclude that they fully characterize the projected Llull matrices, and that the operator $(v_{xy})\mapsto(\psc_{xy})$ really deserves being called a projection:

\smallskip
\begin{theorem}\hskip.5em
\label{st:projection}
The operator $P:\Omega\ni(v_{xy})\mapsto(\psc_{xy})\in\Omega$ is idempotent, \ie $P^2=P$. Its image $P\Omega$ consists exactly of the Llull matrices $(v_{xy})$ that satisfy $(\ref{eq:vxyinequality}\text{--}\ref{eq:ttm})$ for some total order $\xi$.
\end{theorem}

\section{The rank-like rates}

Let us recall that the rank-like rates $\rlr x$ are determined
from the projected scores by formula~(\ref{eq:rrates}).
The properties of the projected scores obtained in Theorem~\ref{st:propertiesOfProjection} imply the following facts, which admit the same proof as in~\cite{crc}
(the only change is that Theorem~\ref{st:propertiesOfProjection} of the present article must be invoked instead of \cite[Thm.~6.3]{crc}):

\smallskip
\newcommand\bla{Same proof as in \cite[Lem.~7.1]{crc}}
\begin{lemma}[\bla]\hskip.5em
\label{st:obsRosa}

\iim{a} If $x\rxi y$ in an admissible order $\xi$,
then $\rlr{x}\le\rlr{y}$.

\iim{b} $\rlr{x}=\rlr{y}$ \,\ifoi\, $\psc_{xy}=\psc_{yx}$.

\iim{c} $\rlr{x}\le\rlr{y}$ \,implies\ \,
the inequalities
$(\ref{eq:pvxyinequality})$ and $(\ref{eq:pvinequalities})$.

\iim{d} $\rlr{x}<\rlr{y}$ \,\ifoi\, $\psc_{xy}>\psc_{yx}$.

\iim{e} $\psc_{xy}>\psc_{yx}$ \,implies\ \, $x\rxi y$ in any admissible order $\xi$.
\end{lemma}

\vskip2pt
\renewcommand\bla{Same proof as in \cite[Thm.~7.2]{crc}}
\begin{theorem}[\bla]\hskip.5em
\label{st:RvsNu}
The rank-like rating given by~\textup{(\ref{eq:rrates})}
is related to the indirect comparison relation
$\icr$ in the following way:
\begin{gather}
\label{eq:referee1}
\rlr{x} < \rlr{y} \ensep\Longleftrightarrow\ensep yx\not\in(\hat\icr)^*,
\\[2.5pt]
\label{eq:referee1bis}
\rlr{x} \le \rlr{y} \ensep\Longleftrightarrow\ensep xy\in(\hat\icr)^*,
\end{gather}
where $\hat\icr$ is the codual of $\icr$, namely $\hat\icr = \{xy\mid\isc_{xy} \ge \isc_{yx}\}$.
\end{theorem}

\vskip2pt
\renewcommand\bla{Same proof as in \cite[Cor.~7.3]{crc}}
\begin{corollary}[\bla]\hskip.5em
\label{st:RvsNuCor}

\iim{a}$\rlr{x} < \rlr{y} \,\Rightarrow\, xy\in\icr$.


\iim{b}If $\hat\icr$ is transitive (in particular, if $\icr$ is total),
then $\rlr{x}\!<\!\rlr{y} \kern-.25pt\Leftrightarrow\kern-.25pt xy\!\in\!\icr$.

\iim{c}If $\icr$ contains a set of the form $\xst\times\yst$ with $\xst\cup\yst=\ist$, then\, $\rlr{x} < \rlr{y}$ for any $x\in\xst$ and $y\in\yst$.
\end{corollary}

\medskip
In the complete case, the fact that $\pto_{xy}=\psc_{xy}+\psc_{yx}=1$ implies that $\sum_{x\in\ist} \rlr{x} = N(N+1)/2$. Related to it, 
one has the following general fact:




\begin{lemma}\hskip.5em 
\label{p:mitjana-rangs}
For any $\xst\subseteq\ist$ one has
\begin{equation}
\label{eq:sumarangs}
\sum_{x\in\xst} \rlr{x} \,\ge\, |\xst|\,(|\xst|\cd+1)/2.
\end{equation}
This inequality becomes an equality 
\,when and only when\, the two following conditions are satisfied:
\begin{align}
\label{eq:mitjana-rangs-2}
\pto_{x\bar x}\,&=\,1,\qquad \hbox{for all\, $x,\bar x\in\xst$}
\\[2.5pt]
\label{eq:mitjana-rangs-3}
\psc_{xy}\,&=\,1,\qquad \hbox{for all\, $x\in\xst$ and $y\not\in\xst$.}
\end{align}
\end{lemma}

\begin{proof}
Starting from formula~(\ref{eq:rrates}), we obtain
\begin{equation*}
\begin{split}
\label{ineq:mitjana-rangs-1}
\sum_{x\in\xst} \rlr{x} \,&=\, N |\xst|
\,-\, \sum_{\latop{\scriptstyle x,\bar x\in\xst}{\scriptstyle \bar x\neq x}} \psc_{x\bar x}
\,-\, \sum_{\latop{\scriptstyle x\in\xst}{\scriptstyle y\not\in\xst}} \psc_{xy} \\[2.5pt]
\,&\ge\, N |\xst| \,-\, |\xst|\,(|\xst|\cd-1)/2 \,-\, |\xst|\,(N\cd-|\xst|) \\[2.5pt]
\,&=\, |\xst|\,(|\xst|\cd+1)/2,
\end{split}
\end{equation*}
where the inequality derives from the following ones: $\pto_{x\bar x} = \psc_{x\bar x}+\psc_{\bar x x}\le 1$\linebreak for $x,\bar x\in\xst$, and $\psc_{xy}\le 1$ for $x\in\xst$ and $y\not\in\xst$.
\end{proof}

\section{Continuity}

We claim that the rank-like rates $\rlr x$ are continuous functions of the binary scores~$v_{xy}$. The main difficulty in proving this statement lies in the admissible order~$\xi$, which plays a central role in the computations. Since $\xi$ varies in a discrete set, its dependence on the data cannot be continuous at~all. Even so, we claim that the final result is still a~continuous function of the data.

In this connection, one can consider as data the normalized Llull matrix~$(v_{xy})$, its domain of variation being the set $\Omega$ introduced in \secpar{1.3}.
\ensep
Alternatively, one can consider as data the relative frequencies of the possible votes, \ie the coefficients $\alpha_k$ mentioned also in \secpar{1.3}.

\smallskip
\begin{theorem}\hskip.5em
\label{st:continuityThm} The~projected scores~$\psc_{xy}$ and
the~rank-like rates~$\rlr{x}$ depend continuously on the Llull
matrix~$(v_{xy})$.
\end{theorem}

\begin{proof}\hskip.5em
Let us begin by considering the dependence of the rank-like rates on the projected scores. This dependence is given by formula~(\ref{eq:rrates}), which is not only continuous but even linear (non-homogeneous).\ensep 

So we are left with the problem of showing that the projection $P:(v_{xy})\mapsto(\psc_{xy})$ is continuous. By arguing as in \cite[Theorem~8.1]{crc}, the problem reduces to showing the continuity of $P_\xi:\Omega_\xi\rightarrow\Omega$ for an arbitrary total order~$\xi$, where $\Omega_\xi$ means the subset of~$\Omega$ which consists of the Llull matrices for which $\xi$ is an~admissible order, and $P_\xi$ means the restriction of $P$ to $\Omega_\xi$.

In order to check that $P_\xi$ is continuous for every total order $\xi$, one has to go over the different mappings whose composition defines $P_\xi$ (see~\secpar{2.1}), namely:\ensep
$(v_{xy})\mapsto(\isc_{xy})\mapsto(\img_{xy})$,\ensep
$(v_{xy})\mapsto(t_{xy})$,\ensep
$(\img_{xy})\mapsto(\ppmg_{xy})$,\ensep
$\Psi:((\ppmg_{xx'}),(t_{xy}))\mapsto(\ppto_{xy})$,\ensep
and finally $((\ppmg_{xx'}),(\ppto_{xx'}))\mapsto(\psc_{xy})$.
\ensep
Except for $\Psi$, all of these mappings involve only the operations of addition, subtraction, multiplication by a constant, maximum and minimum, which are certainly continuous. Concerning the operator $\Psi$, let us recall that its output is the orthogonal projection of $(t_{xy})$ onto a certain convex set determined by $(\ppmg_{xx'})$; a general result of continuity for such an operation can be found in~\cite{da}.
\end{proof}


\smallskip
\begin{corollary}\hskip.5em
\label{st:continuityCor} The rank-like rates  depend continuously on
the relative frequency of each possible content of an~individual
vote.
\end{corollary}

\begin{proof}\hskip.5em
It suffices to notice that the Llull matrix $(v_{xy})$ is simply the center of gravity of the distribution specified by these relative frequencies (formula~(\ref{eq:cog}) of~\secpar{1.3}).
\end{proof}

\section{Decomposition}

The decomposition property~\llrd\ 
is concerned with having a partition of~$\ist$ into two non-empty sets $\xst$ and $\yst$ such that each member of $\xst$ is unanimously preferred to any member of $\yst$, that is:
\begin{equation}
\label{eq:v1}
v_{xy} \,=\, 1
\quad \hbox{(and therefore $v_{yx}=0$)}
\quad \hbox{whenever\, $xy\in\xst\times\yst$}.
\end{equation}
More specifically, property~\llrd, which was proved in \cite[\secpar{9}]{crc}, ensures that in the complete case
such a situation is characterized by any of the following equalities:
\begin{alignat}{2}
\rlr{x} \,&=\, \rlrbis{x},\qquad &&\hbox{for all $x\in\xst$},
\label{eq:condrfx}\\[3.5pt]
\rlr{y} \,&=\, \rlrbis{y} \,+\, |X|,\qquad &&\hbox{for all $y\in\yst$},
\label{eq:condrfy}\\
\sum_{x\in\xst}\rlr{x} \,&=\, |\xst|\,(|\xst|+1)/2 &&\null,
\label{eq:condrfs}
\end{alignat}
where $\rlrbis{x}$ and $\rlrbis{y}$ denote the rank-like rates which are determined respectively from the submatrices associated with $\xst$ and $\yst$.

\smallskip
In this section we will see that some of these implications are still valid under certain assumptions that allow for incompleteness.
In particular, Theorem~\ref{st:mitjana} below entails that under the assumption of transitive individual preferences an option gets a rank-like rate exactly equal to~$1$ \ifoi it is unanimously preferred to any other.

\renewcommand\upla{\vskip-2pt}

\upla
\medskip
\begin{lemma}\hskip.5em
\label{lema1-vxy}
Given a partition $\ist=\xst\cup\yst$ into two disjoint nonempty sets, one has the following implications:
\begin{equation}
\left.
\begin{array}{c} v_{xy}=1 \\[2.5pt]
\forall\,xy\in\xst\times\yst
\end{array}
\right\}
\ \Longrightarrow\
\left\{
\begin{array}{c} \img_{xy}\cd=1 \\[2.5pt]
\forall\,xy\in\xst\times\yst
\end{array}
\right\}
\ \Longleftrightarrow\
\left\{
\begin{array}{c}
v_{xy}^\pi=1 \\[2.5pt]
\forall\,xy\in\xst\times\yst
\end{array}
\right. \label{eq:lema1-vxy-general}
\end{equation}
If the individual preferences are transitive,
\,then the converse of the first implication holds too.
\end{lemma}

\begin{proof}\hskip.5em
Here we will only prove the converse of the first implication in the case of transitive individual preferences. All the other statements are proved by the arguments given in  \cite[Lemma~9.1]{crc}.

Assume that $\img_{xy}=1$ for all $xy\in\xst\times\yst$.
Since $\img_{xy}=\isc_{xy}-\isc_{yx}$ and both terms of
this difference belong to $[0,1]$, the only possibility is
$\isc_{xy}=1$ (and $\isc_{yx}=0$). This implies the existence of
a path $x_0x_1\dots x_n$ from $x$ to $y$ such that $v_{x_ix_{i+1}} =
1$ for all $i$. But this means that all of the votes include each of
the pairs $x_ix_{i+1}$ of this path. So, if they are transitive, all~of them include also the pair $xy$, \ie $v_{xy} = 1$.
\end{proof}

\upla
\medskip
\begin{theorem}\hskip.5em
\label{st:mitjana}
Condition $(\ref{eq:v1})$ implies $(\ref{eq:condrfs})$.
If the individual preferences are transitive, 
\,then the converse implication holds too.
\end{theorem}

\begin{proof}\hskip.5em
On account of Lemmas~\ref{lema1-vxy} and \ref{p:mitjana-rangs}, it suffices to see that (\ref{eq:mitjana-rangs-3}) implies (\ref{eq:mitjana-rangs-2}).
So, let us assume that $\psc_{xy}=1$ for any $x\in\xst$ and $y\in\yst$. By~Lemma~\ref{st:obsRosa}.(e), any admissible order $\xi$ includes all pairs $xy$ with $x\in\xst$ and $y\in\yst$. Let $\ell$ be the last element of $\xst$ in a fixed admissible order. Then $\psc_{\ell\ell'}=1$ and therefore $\pto_{\ell\ell'}=1$. On account of the pattern of growth of the projected turnouts, that is inequalities~(\ref{eq:ptinequalities}), this implies that $\pto_{x\bar x}=1$ for any $x,\bar x\in\xst$. 
\end{proof}

\upla
\smallskip
\begin{corollary}\hskip.5em
If the individual preferences are transitive,
\,then $\rlr{x}=1$ \ifoi $v_{xy}=1$ for every $y\ne x$.
\end{corollary}

\upla
\medskip
\begin{theorem}\hskip.5em
\label{st:condrfPro} If the individual votes are rankings (possibly
truncated or with ties),\, then $(\ref{eq:v1})$ implies $(\ref{eq:condrfx})$.
\end{theorem}

\begin{proof}\hskip.5em
As in~(\ref{eq:condrfx}) we will continue using a tilde to distinguish between hom\-olo\-gous objects associated respectively with the whole matrix and with the submatrix associated with $\xst$.

First of all we will show that
\begin{equation}
\label{eq:xa}
t_{xy} \,=\, 1,\qquad\hbox{for all $xy\in\xst\times\ist$.}
\end{equation}
In fact, the rules that we
are using for translating votes into binary preferences ---namely,
rules~(a--d) of \secpar{1.3}--- entail the following
implications:\ensep (i)~$v_{xy}=1$ for some $y\in\ist$ implies that
$x$~is explicitly mentioned in all of the ranking votes;\ensep
and\,~(ii)~$x$~being explicitly mentioned in all of the ranking
votes implies that $t_{xy}=1$ for any $y\in\ist.$ 

In~particular, (\ref{eq:xa}) 
ensures that the Llull matrix restricted to $\xst$ is complete,
from which it follows that
\begin{equation*}
\label{eq:totals-restringits-x}
\ptobis_{x\bar x} \,=\, 1,\qquad\hbox{for all $x,\bar x\in\xst$.}
\end{equation*}

Concerning the non-restricted matrix, we know, by Lemma~\ref{lema1-vxy}, that 
condition~(\ref{eq:v1}) implies $\psc_{xy}=1$ and therefore $\pto_{xy}=1$ for all $xy\in\xst\times\yst$. As in the proof of Theorem~\ref{st:mitjana}, this implies that
\begin{equation*}
 \label{eq:totals-x}\pto_{x\bar x} \,=\, 1,\qquad\hbox{for all $x,\bar x\in\xst$.}
\end{equation*}

On the other hand, Lemma~9.2 of \cite{crc}, whose proof is valid in the general case,
ensures that $\ppmgbis_{x\bar x}=\ppmg_{x\bar x}$ for all $x,\bar x\in\xst$.
\ensep
Altogether, we get $\pscbis_{x\bar x}=\psc_{x\bar x}$ for all $x,\bar x\in\xst$.
\ensep
Finally, (\ref{eq:condrfx}) is a direct consequence of this equality together with the above-remarked fact that $\psc_{xy}=1$ for all $xy\in\xst\times\yst$.
\end{proof}

\section{Other properties}

In this section we collect several other properties whose proof given in~\cite{crc} remains valid in the general case. The only caveat to bear in mind is that here they ultimately rely on Theorem~\ref{st:propertiesOfProjection} of the present article instead of \cite[Thm.~6.3]{crc}. One of these properties is here complemented by an additional result that was not mentioned in~\cite{crc}. 

\bigskip
The first of these properties is the Condorcet-Smith principle~\llmp:

\renewcommand\bla{Same proof as in \cite[Thm.~10.1]{crc}}
\begin{theorem}[\bla]
\label{st:majPrinciple} Both the indirect majority
relation
$\icr$ and the ranking 
determined by the rank-like rates comply with the Condorcet-Smith principle:\ensep
If $\ist$ is partitioned into two sets $\xst$ and $\yst$ with the
property that $v_{xy} > 1/2$ for any $x\in\xst$ and $y\in\yst$,
\ensep
then one has also $xy\in\icr$ and $\rlr{x}<\rlr{y}$
for any such $x$ and $y$.
\end{theorem}

\bigskip
The next results are concerned with clone consistency.
In this connection we make use of the notion of autonomous sets.
\ensep
A subset $\cst\sbseteq\ist$ is said to be {\df autonomous} for a relation~$\rho$ when
each element from outside~$\cst$ relates to all elements of~$\cst$ in the same way;
in other words, when, for any~$x\not\in\cst$, having
$ax\in\rho$ for some $a\in\cst$ implies $bx\in\rho$ for any
$b\in\cst$, and similarly, having $xa\in\rho$ for some $a\in\cst$
implies $xb\in\rho$ for any $b\in\cst$.
\ensep
More generally, a subset $\cst\sbseteq\ist$ will be said to be autonomous for a valued relation $(v_{xy})$ when
the equalities $v_{ax} = v_{bx}$ and $v_{xa} = v_{xb}$ hold whenever
$a,b\in\cst$ and $x\not\in\cst$.
\ensep
For more details about the notion of autonomous set and the property of clone consistency we refer the reader to \cite[\secpar{11}]{crc}.

\smallskip
\renewcommand\bla{Same proof as in \cite[Thm.~11.5 and 11.7]{crc}}
\begin{theorem}[\bla] 
\label{st:clons4}
Assume that $\cst\sbset\ist$ is autonomous for the 
Llull matrix $(v_{xy})$. Then $\cst$ is autonomous for
the indirect comparison relation
$\icr$ as well as for the ranking 
determined by the rank-like rates, \ie for the
relation $\hat\rlrating = \{xy\in\tie\mid \rlr{x}\le\rlr{y}\}$.
\ensep
Besides, contracting $\cst$ to a single option in the Llull matrix
has no other effect in~$\icr$ and $\hat\rlrating$ than getting the same contraction. 
\end{theorem}

\vskip-6pt
\remark In contrast to \cite[Thm.~11.6]{crc}, in the incomplete case $\cst$ is not ensured to be autonomous for the projected scores $(\psc_{xy})$. In fact, although the intermediate projected margins $(\ppmg_{xy})$ do have this property, the intermediate projected turnouts $(\ppto_{xy})$ can do away with it.
It would certainly be interesting to find an alternative to step 4 (quadratic minimization with constraints) so that the projected scores $(\psc_{xy})$ keep the autonomous sets of the original Llull matrix $(v_{xy})$ (in addition to the present properties).

\vskip4pt
\medskip
Together with the facts stated in the preceding theorem, 
one could expect that the restriction of the final ranking to $\cst$
should coincide with the ``local'' result that one obtains
when starting from the restriction of the original Llull matrix to $\cst\times\cst$.
A requirement of this sort (though concerning only the winner) is included for instance
in the property of ``composition consistency'' considered in~\cite{lll}.
Our method does not entirely satisfy it.
This is due to the fact that $\cst$ being autonomous does not prevent
the indirect score $\isc_{ab}$ for $a,b\in\cst$ to come from outside $\cst$.
However, the next result still ensures that the differences between the local ranking and the global one are limited to some strict preferences of the former being replaced by ties in the latter.

\smallskip
\begin{theorem} 
\label{st:clons5}
Assume that $\cst\sbset\ist$ is autonomous for the Llull matrix~$(v_{xy})$.
Let $\rlrbis{x}$ denote the rank-like rates that are obtained on $\cst$
when starting from the restriction of $(v_{xy})$ to $xy\in\cst\times\cst$.
The ranking that these rates determine on $\cst$ is related to that determined by the global rates $\rlr{x}$ in the following way:
\begin{equation}
\label{eq:octubre2011a}
\rlrbis{a} \le \rlrbis{b} \,\Longrightarrow\, \rlr{a} \le \rlr{b},\qquad\hbox{whenever $a,b\in\cst$.}
\end{equation}
\end{theorem}

\begin{proof}\hskip.5em
We will systematically use a tilde to denote the objects that are obtained when starting from the restriction of $(v_{xy})$ to $xy\in\cst\times\cst$. Let $a,b$ be two arbitrary elements of $\cst$. The implication (\ref{eq:octubre2011a}) will be a consequence of the following one:
\begin{equation}
\label{eq:octubre2011b}
\iscbis_{ab} \ge \iscbis_{ba} \,\Longrightarrow\, \isc_{ab} \ge \isc_{ba},\qquad\hbox{whenever $a,b\in\cst$.}
\end{equation}
In fact, this is saying that $\hat\icrbis\sbseteq\hat\icr$, which obviously entails $(\hat\icrbis)^*\sbseteq(\hat\icr)^*$, and therefore gives (\ref{eq:octubre2011a}) by virtue of 
Theorem~\ref{st:RvsNu}.

In order to prove (\ref{eq:octubre2011b}) we begin by noticing that the definition of the indirect scores ensures the inequality $\isc_{ab}\ge\iscbis_{ab}$, as well as the analogous one for $ba$ instead of $ab$. This immediately settles (\ref{eq:octubre2011b}) in the case $\isc_{ba}=\iscbis_{ba}$. In the case $\isc_{ba}>\iscbis_{ba}$ one can argue as follows. Having this strict inequality means that
the maximum that defines $\isc_{ba}$ is achieved by a path $\beta$ from $b$ to $a$ that contains some $x\in\ist\setminus\cst$. Let $\beta_1$ and $\beta_2$ denote the segments of this path that go respectively from $b$ to $x$ and from $x$ to $a$.
The desired result is then obtained in the following way, where the second step makes use of the fact that 
$\cst$ is autonomous for the indirect scores \cite[Prop.~11.3]{crc}:
\begin{equation}
\isc_{ab} \,\ge\, \min(\isc_{ax},\isc_{xb})
 \,=\, \min(\isc_{bx},\isc_{xa})
 \,\ge\, \min(v_{\beta_1},v_{\beta_2}) \,=\, v_\beta \,=\, \isc_{ba}.
\end{equation}
\end{proof}

\medskip
Finally, we have the following (rather weak) property of monotonicity:

\smallskip
\renewcommand\bla{Same proof as in \cite[Thm.~12.1 and Cor.~12.2]{crc}}
\begin{theorem}[\bla]\hskip.5em
\label{st:mono} Assume that $(v_{xy})$ and $(\vbis_{xy})$ are
related to each other in the following way:
\begin{equation}
\label{eq:mona}
\vbis_{ay} \ge v_{ay},\quad \vbis_{xa} \le v_{xa},\quad \vbis_{xy} = v_{xy},\qquad \forall x,y\neq a.
\end{equation}
In this case, the following properties are satisfied for any
$x,y\neq a$:
\begin{gather}
\iscbis_{ay} \ge \isc_{ay},\qquad \iscbis_{xa} \le \isc_{xa},
\label{eq:mono1}
\\[2.5pt]
\rlr{a} < \rlr{y} \,\Longrightarrow\, \rlrbis{a} \le \rlrbis{y},
\label{eq:mono3}
\\[2.5pt]
\label{eq:smona} (\rlr{a} < \rlr{y},\ \forall y\neq a)\ \,\Longrightarrow\,\
(\rlrbis{a} < \rlrbis{y},\ \forall y\neq a).
\end{gather}
\end{theorem}


\section{Approval voting}\label{sec-approval}

In approval voting, each voter is asked for a list of approved
options, without any expression of preference between them, and each
option~$x$ is then rated by the number of approvals for
it~\cite[\secpar{1 and 2}]{br}. In the following we will refer to this number as the
\dfc{approval score} of~$x$, and its value relative to the total number of votes $V$ will be denoted by $\av{x}$.

From the point of view of paired comparisons, an individual vote of
approval type can be viewed as a truncated ranking where all the
options that appear in it are tied. In~this section, we will see
that the margin-based variant
orders the options exactly in the same way as the approval scores.
In~other words, the method of approval voting agrees with ours in a qualitative way under
interpretation~(d$'$) of~\secpar{1.3},
\ie under the interpretation that the non-approved options of each
individual vote are tied.
This interpretation is in congruence with the
hypothesis of dichotomous preferences,
under which approval voting has especially good properties~\cite{vorsatz}.

Having said that, the preliminary result~\ref{st:av1} 
will hold not only under interpretation~(d$'$) but also under
interpretation~(d), \ie that there is no information about the
preference of the voter between two non-approved options, and also
under the analogous interpretation that there is no information
about his preference between two approved options.\ensep
Interpretation~(d$'$) does not play an essential role until
Theorem~\ref{st:av4}, where we use the fact that it always brings
the problem into the complete case.

\smallskip
In the following we use the following notation:
\begin{alignat}{2}
&\mu(v) \,=\, \{xy\mid v_{xy} > v_{yx}\},\qquad 
&&\hat\mu(v) \,=\, \{xy\mid v_{xy} \ge v_{yx}\},
\\[2.5pt]
&\mu(\av{}) \,=\, \{xy\mid \av{x} > \av{y}\},\qquad
&&\hat\mu(\av{}) \,=\, \{xy\mid \av{x} \ge \av{y}\}.
\end{alignat}
Notice that $\mu(\isc)$ is the indirect comparison relation $\icr$.

\medskip
\begin{proposition}\hskip.5em
\label{st:av1} In the approval voting situation, the following
equality holds:
\begin{equation}
v_{xy} - v_{yx} \,=\, \av{x} - \av{y}. \label{eq:mgav}
\end{equation}
As a consequence, $\mu(v)=\mu(\av{})$.
\end{proposition}
\begin{proof}\hskip.5em
Obviously, the possible ballots are in one-to-one correspondence
with the subsets $\xst$ of $\ist$. In the following, $v_\xst$
denotes the relative number of votes that approved exactly the set
$\xst$. With this notation it is obvious that
\begin{equation}
\av{x} \,=\, \sum_{\xst\ni\,x} v_{\xst} \,=\,
\sum_{\latop{\scriptstyle \xst\ni\,x}{\scriptstyle \xst\not\ni\,y}}
v_{\xst} + \sum_{\latop{\scriptstyle \xst\ni\,x}{\scriptstyle
\xst\ni\,y}} v_{\xst},
\end{equation}
for any $y\in\ist$.
On the other hand, one has
\begin{equation}
v_{xy} \,=\, \sum_{\latop{\scriptstyle \xst\ni\,x}{\scriptstyle
\xst\not\ni\,y}} v_{\xst} \,+\, \bigg[\,{\textstyle\frac12}
\,\sum_{\latop{\scriptstyle \xst\ni\,x}{\scriptstyle \xst\ni\,y}}
v_{\xst} + {\textstyle\frac12} \,\sum_{\latop{\scriptstyle
\xst\not\ni\,x}{\scriptstyle \xst\not\ni\,y}} v_{\xst}\,\bigg],
\label{eq:vxyav}
\end{equation}
where each of the terms in brackets is present or not depending on which interpretation is used.
Anyway, the preceding expressions, together with the analogous ones
where $x$ and $y$ are interchanged with each other, result in the
equality~(\ref{eq:mgav}) independently of those alternative
interpretations.
\end{proof}

\medskip
\begin{lemma}\label{l:creix-up-right} Assume that there exists a total order $\xi$ such that the scores~$(v_{xy})$ satisfy \,$\mu(v)\subseteq \xi$\, together with the condition
\begin{equation}\label{eq:creix-up-right}
v_{xz} \ge v_{yz},\quad v_{zx}\le v_{zy},\qquad \hbox{whenever $x\rxi y$ and $z\not\in\{x,y\}$.}
\end{equation}
Then one has also $\mu(\isc)\subseteq \xi$, \ie $\xi$ is an admissible order.
In the complete case one has the equality \,$\mu(\isc)=\mu(v).$
\end{lemma}

\begin{proof}\hskip.5em
Let us begin by recalling that $\mu(v)\sbseteq\xi$ is equivalent to $\xi\sbseteq\hat\mu(v)$, and similarly for $\mu(\isc)$ instead of $\mu(v)$.
\ensep
The first statement of the lemma will be obtained by showing that under its hypotheses one has
\begin{equation}
\label{eq:vvvv}
\isc_{xy} \,=\, v_{xy} \,\ge\, \isc_{yx} \,\ge\, v_{yx},\qquad\hbox{whenever $x\rxi y$.}
\end{equation}
On account of the definition of $\isc_{xy}$ and the fact that $v_{xy}\le\isc_{xy}$ (and analogously for $yx$), in order to prove~(\ref{eq:vvvv}) it suffices to show that $x\rxi y$ implies
\begin{equation}
\label{eq:vvvvgamma}
v_\gamma \,\le\, v_{xy},\quad v_\eta \,\le\, v_{xy},
\end{equation}
for any path $\gamma=x_0x_1\dots x_n$ from $x_0=x$ to $x_n=y$,
and for any path $\eta=y_0y_1\dots y_n$ from $y_0=y$ to $y_n=x$.
Without loss of generality, in the following we will assume $n>1$ and we will let $\eta$ be the reverse of $\gamma$, \ie $y_i = x_{n-i}$.
\ensep
Let us assume that $x\rxi y$. In~order to prove (\ref{eq:vvvvgamma}) we will distinguish three cases, namely:\ensep
(i)~$x\rxieq x_i\rxieq y$ for all~$i$;\ensep
(ii)~$x_i\rxi x$ for some~$i$;\ensep
(iii)~$y\rxi x_i$ for some~$i$.

Case~(i)\,: It suffices to notice that the definition of the score of a path and the assumptions of the lemma allow to write $v_\gamma \le v_{xx_1} \le v_{xy}$ and also $v_\eta \le v_{x_1x} \le v_{xx_1} \le v_{xy}$. For future reference, let us notice also that in the complete case with $v_{xy}>1/2$ one can write $v_\eta \le v_{x_1x} \le 1/2 < v_{xy}$, so $v_\eta$ is then strictly less than $v_{xy}$.

Case~(ii)\,: Let $i$ be the first time that one has $x_i\rxi x$, and let $j$ be the first time after $i$ that one has $x\rxieq x_j$. Obviously, $0<i<j\le n$. By~construction, we have $x\rxieq x_{i-1}$ and $x_i\rxi x\rxi y$, which entail that $v_\gamma \le v_{x_{i-1}x_i} \le v_{xx_i} \le v_{xy}$. On the other hand, we have $x_{j-1}\rxi x\rxi y$ and $x\rxieq x_j$, which entail $v_\eta \le v_{x_jx_{j-1}} \le v_{xx_{j-1}} \le v_{xy}$. Similarly as above, in the complete case with $v_{xy}>1/2$ we get the strict inequality $v_\eta<v_{xy}$.

Case~(iii)\, is analogous to case~(ii).

\smallskip
The statement about the complete case will be proved if we show the following implications
\begin{gather}
\label{eq:vvvva}
v_{xy} \,>\, v_{yx} \ensep\Longrightarrow\ensep \isc_{xy} \,>\, \isc_{yx},
\\[2.5pt]
\label{eq:vvvvb}
v_{xy} \,=\, v_{yx} \ensep\Longrightarrow\ensep \isc_{xy} \,=\, \isc_{yx}.
\end{gather}
In order to obtain (\ref{eq:vvvva}) it suffices to use the already remarked fact that in the present circumstances the second inequality of (\ref{eq:vvvvgamma}) is strict, which implies that one has also a strict inequality in the middle of (\ref{eq:vvvv}).
\ensep
Finally, in order to obtain (\ref{eq:vvvvb}) it suffices to use (\ref{eq:vvvv}) with $v_{xy}=v_{yx}=1/2$.
\end{proof}

\medskip
\begin{theorem}\hskip.5em
\label{st:av4} In the approval voting situation, the margin-based
variant results in a full qualitative compatibility between the rank-like rates $\rlr x$ and the approval scores $\av{x}$, in the sense that \,$\rlr x <
\rlr y \,\Longleftrightarrow\, \av{x} > \av{y}$.
\end{theorem}

\begin{proof}\hskip.5em
Recall that the margin-based variant amounts to using
interpretation~(d$'$), which always brings the problem into the
complete case (when the terms in brackets are included,
equation~(\ref{eq:vxyav}) has indeed the property that
$v_{xy}+v_{yx}=1$).

Let $\xi$ be any total ordering of the elements of~$\ist$ by non-increasing values of~$\av{x}$. In other words, $\xi$ is any total order contained in $\hat\mu(\av{})$, which is equivalent to say, any total order containing $\mu(\av{})$. We claim that we are under the hypothesis of the preceding lemma. This is a consequence of Proposition~\ref{st:av1}.
\ensep
In fact, on the one hand it immediately gives $\mu(v)\sbseteq\xi$.
\ensep
On the other hand, it allows to deal with the margins $m_{xy} = v_{xy}-v_{yx}$ in the following way:
\begin{alignat*}{5}
&m_{xz} &&\,=\, \av{x}-\av{z} &&\,=\, (\av{x}-\av{y})+(\av{y}-\av{z})&&\,=\, m_{xy}+m_{yz} &&\,\ge\, m_{yz},
\\[2.5pt]
&m_{zx} &&\,=\, \av{z}-\av{x} &&\,=\, (\av{z}-\av{y})+(\av{y}-\av{x}) &&\,=\, m_{zy}+m_{yx} &&\,\le\, m_{zy},
\end{alignat*}
where we have assumed $x\rxi y$ and we have used that $m_{xy}=\av{x}-\av{y}\ge 0$ and $m_{yx}=-m_{xy}\le 0$. The obtained inequalities are certainly equivalent to those of~(\ref{eq:creix-up-right}).

By virtue of Lemma~\ref{l:creix-up-right}, we are therefore ensured that $\icr := \mu(\isc) = \mu(v)$. Since we know that $\mu(v)=\mu(\av{})$, we can also write $\icr = \mu(\av{})$ and $\hat\icr = \hat\mu(\av{})$, which guarantees that $\hat\icr$ is transitive. This allows to apply  Corollary~\ref{st:RvsNuCor}.(b) to~arrive at the conclusion that
$$
\rlr x<\rlr y
\,\Longleftrightarrow\, xy\in \icr \,\Longleftrightarrow \,\av{x}>\av{y}.
$$
\end{proof}

\vskip-15mm\null 

\end{document}